\documentclass{article}
\usepackage{hyperref}
\usepackage{amsthm,amsmath}
\usepackage{amssymb}
\usepackage{bm}
\usepackage{bbm}
\usepackage{url}
\usepackage{comment}
\usepackage{url}
\usepackage{multirow}
\usepackage{graphicx}
\usepackage{color}
\usepackage{tikz}
\usepackage{graphicx}
\usepackage{caption}
\usepackage{float}
\usepackage{mathtools}
\usepackage{pifont}
\usepackage{amsmath, enumerate, algorithm,algorithmic,caption,subcaption,graphicx,color}
\usepackage{hyperref}
\usepackage{fullpage}

\allowdisplaybreaks

\title{Semi-proximal Mirror-Prox \\for Nonsmooth Composite Minimization
\thanks{The authors would like to thank Anatoli Juditsky and Arkadi Nemirovski for fruitful discussions. 
This work was supported by the NSF Grant CMMI-1232623, the LabEx Persyval-Lab (ANR-11-LABX-0025),
   the project Titan (CNRS-Mastodons),
the project Macaron (ANR-14-CE23-0003-01),
the MSR-Inria joint centre,
and the Moore-Sloan Data Science Environment at NYU.}}
\author{Niao He \\
\texttt{nhe6@gatech.edu} \\
GeorgiaTech
\and
Zaid Harchaoui \\
\texttt{zaid.harchaoui@inria.fr} \\
NYU, Inria
}

\def\cI{{\cal I}}
\def\SadVal{\hbox{\rm SadVal}}
\def\epsilonsad{{\epsilon_{\hbox{\tiny\rm Sad}}}}
\def\mybox\blacksquare
\def\Bbb#1{{\bf #1}}

\def\fnote#1{\footnote}
\def\blacksquare{\hbox{\vrule width 4pt height 4pt depth 0pt}}

\def\beq{\begin{equation}}
\def\eeq{\end{equation}}

\newtheorem{theorem}{Theorem}[section]

\newtheorem{corollary}{Corollary}[section]

\newtheorem{proposition}{Proposition}[section]

\def\Argmin{\mathop{\rm Argmin}}

\def\Argmin{\mathop{\rm Argmin}}

\def\hat{\widehat}

\def\Bbb#1{{\bf #1}}

\def\fnote#1{\footnote}
\def\blacksquare{\hbox{\vrule width 4pt height 4pt depth 0pt}}

\def\beq{\begin{equation}}
\def\eeq{\end{equation}}

\def\ba{\begin{array}}
\def\ea{\end{array}}
\def\beann{\begin{eqnarray*}}
\def\eeann{\end{eqnarray*}}
\def\bea{\begin{eqnarray}}
\def\eea{\end{eqnarray}}

\def\Argmin{\mathop{\rm Argmin}}

\def\Bbb#1{{\bf #1}}

\def\fnote#1{\footnote}
\def\blacksquare{\hbox{\vrule width 4pt height 4pt depth 0pt}}

\def\Argmin{\mathop{\rm Argmin}}

\def\cB{{\cal B}}

\def\cN{{\cal N}}

\def\epsilonvi{{\epsilon_{\hbox{\scriptsize\rm VI}}}}

\def\Argmin{\mathop{\hbox{\rm  Argmin}}}

\def\Bbb#1{{\bf #1}}
\def\beq{\begin{equation}}

\def\blacksquare{\hbox{\vrule width 4pt height 4pt depth 0pt}}
\def\bR{{\mathbb{R}}}

\def\eeq{\end{equation}}

\def\fnote#1{\footnote}

\def\mypict3{\epsfxsize=220pt\epsfysize=80pt\epsfbox}

\def\Opt{{\hbox{\rm  Opt}}}

\def\Res{{\hbox{\rm  Res}}}

\def\cA{{{\cal A}}}

\def\bR{{\mathbf{R}}}
\def\Res{{\hbox{\rm Res}}}
\def\nuc{{\hbox{\scriptsize\rm nuc}}} 
\def\Prox{\hbox{\rm Prox}}

\newcommand{\be}{\begin{eqnarray}}
\newcommand{\ee}[1]{\label{#1}\end{eqnarray}}

\newcommand{\ese}{\end{eqnarray*}}
\newcommand{\bse}{\begin{eqnarray*}}
\newcommand{\rf}[1]{(\ref{#1})}
\newcommand{\half}{ \mbox{\small$\frac{1}{2}$}}

\def\LMO{\hbox{\rm LMO}}

\newcommand{\spmp}{Semi-Proximal Mirror-Prox~}
\newcommand{\Spmp}{Semi-Proximal Mirror-Prox~}
\newcommand{\spmpnospace}{Semi-Proximal Mirror-Prox}

\newcommand{\semivi}{Semi-VI~}

\begin{document}
\maketitle

\begin{abstract}
We propose a new first-order optimisation algorithm to solve high-dimensional non-smooth composite minimisation problems. Typical
examples of such problems have an objective that decomposes into a non-smooth empirical risk part and a non-smooth regularisation penalty.
The proposed algorithm, called Semi-Proximal Mirror-Prox, leverages the Fenchel-type representation of one part of the objective while handling the other part 
of the objective via linear minimization over the domain. The algorithm stands in contrast with 
more classical proximal gradient algorithms with smoothing, 
which require the computation of proximal operators at each iteration and can therefore be impractical for high-dimensional problems. 
We establish the theoretical convergence rate of Semi-Proximal Mirror-Prox, which exhibits the optimal complexity bounds, 
i.e. $O(1/\epsilon^2)$,
for the number of calls to linear minimization oracle. We present promising experimental results showing the interest of the approach
in comparison to competing methods.
\end{abstract}

\section{Introduction}
A wide range of machine learning and signal processing problems 
can be formulated as the minimization of a composite objective:
\begin{equation}\label{nonsmoothproblem0}
\min_{x\in X} \: F(x) := f(x)+ \|\cB x\|
\end{equation}
where $X$ is closed and convex, $f$ is convex and can be either smooth, 
or nonsmooth yet enjoys a particular structure. The term
$\|\cB x\|$ defines a regularization penalty through a norm 
$\|\cdot\|$, and $x\mapsto\cB x$ a linear mapping on a closed convex set $X$.
The function $f$ usually corresponds to an empirical risk, 
that is an empirical average of a possibly non-smooth loss function
 evaluated on a set of data-points,
while $x$ encodes the learning parameters. 
All in all, the objective $F$
has a doubly non-smooth structure. 

In many situations, the 
objective function $F$ of interest
enjoys a favorable structure, namely a so-called 
\emph{Fenchel-type representation}~\cite{cox2013dual,CoMP13,juditsky13variational}:
\begin{equation}\label{sadpoint}
f(x)=\max_{z\in Z}\; \left\{\langle x, Az\rangle -\psi(z)\right\}
\end{equation}
where $Z$ is convex compact subset of a Euclidean space, 
and $\psi(\cdot)$ is a convex function. Sec.~\ref{sec:exps} will give several examples
of such situations. 
Fenchel-type representations
can then be leveraged to use first-order optimisation algorithms. 


A simple first option to minimise $F$ is using the so-called Nesterov smoothing technique~\cite{nesterov2005smooth} along with a proximal gradient algorithm~\cite{Parikh13}, assuming that the proximal operator associated with $X$ is computationally tractable and cheap to compute. However, this is certainly  not the case when considering problems with norms acting in the spectral domain of high-dimensional matrices, such as the matrix nuclear-norm~\cite{J13} and structured extensions thereof~\cite{Chen12,Bach12}. In the latter situation, another option is to use a smoothing technique now with a conditional gradient or Frank-Wolfe algorithm to minimize $F$, assuming that a \emph{a linear minimization oracle} associated 
with $X$ is cheaper to compute than the proximal operator~\cite{cox2013dual,lan2013complexity,pierucci2014smoothing}. 
Neither option takes advantage of the composite structure of the objective (\ref{nonsmoothproblem0}) or handles the case when the linear mapping $\cB$ is nontrivial. 

\paragraph{Contributions}
Our goal in this paper is to propose a new first-order optimization algorithm, called \spmp,
designed to solve the difficult non-smooth composite optimisation problem~(\ref{nonsmoothproblem0}),
which does not require the exact computation of proximal operators. Instead, the \spmp~relies upon
i) Fenchel-type representability of $f$; ii) Linear minimization oracle associated with $\|\cdot\|$ in the domain $X$.
While the Fenchel-type representability of $f$ allows to cure the non-smoothness of $f$, 
the linear minimisation over the domain $X$ allows to tackle the non-smooth regularisation penalty $\|\cdot\|$. 
We establish the theoretical convergence rate of \spmpnospace, which exhibits the \emph{optimal complexity bounds}, 
i.e. $O(1/\epsilon^2)$, for the number of calls to linear minimization oracle. Furthermore, \spmp generalizes
previously proposed approaches and improves upon them in special cases:
\begin{enumerate}
\item Case $\cB \equiv 0$: \spmp does not require assumptions on “favorable geometry” of dual domains $Z$ or simplicity of $\psi(\cdot)$ in (\ref{sadpoint}). 
\item Case $\cB = \mathbb{I}$: \spmp is competitive with previously proposed approaches~\cite{Lan14, pierucci2014smoothing} based on smoothing techniques.  
\item Case of non-trivial $\cB$: \spmp is the first proximal-free or conditional-gradient-type optimization algorithm for (\ref{nonsmoothproblem0}). 
\end{enumerate}

\paragraph{Related work}
The \spmp algorithm belongs the family of conditional gradient algorithms, whose most basic instance is the Frank-Wolfe algorithm for constrained smooth optimization using a linear minimization oracle; see~\cite{J13,b15,bertsekas:2015}. Recently, in~\cite{cox2013dual, juditsky13variational}, the authors consider constrained non-smooth optimisation when the domain $Z$ has a ``favorable geometry'', \textit{i.e.} the domain is amenable to linear minimisation (favorable geometry), and 
establish a complexity bound with $O(1/\epsilon^2)$ calls to the linear minimization oracle. 
Recently, in~\cite{Lan14}, a method called conditional gradient sliding is proposed to solve similar problems, using a smoothing technique, with a complexity bound in $O(1/\epsilon^2)$ for the calls
to the linear minimization oracle (LMO) and additionally a $O(1/\epsilon)$ bound for the linear operator evaluations. Actually, this $O(1/\epsilon^2)$ bound for the LMO complexity can be shown to be indeed \emph{optimal} for conditional-gradient-type or LMO-based algorithms, 
when solving general\footnote{Related research extended such approaches to stochastic or online settings~\cite{hazan2012projection, garber2013linearly, Lan14}; such settings are beyond the scope of this work.} non-smooth convex problems~\cite{lan2013complexity}. 

However, these previous approaches
are appropriate for objective with a non-composite structure. When applied to our problem (\ref{nonsmoothproblem0}), the smoothing would be applied 
to the objective taken as a whole, ignoring its composite structure. 
Conditional-gradient-type algorithms were recently proposed for composite objectives~\cite{dudik:harchaoui:malik2012,harchaoui2013conditional,ZYS12,pierucci2014smoothing, mu2014scalable}, but cannot be applied for our problem. In~\cite{harchaoui2013conditional},  $f$ is smooth and $\cB$ is identity matrix, whereas in~\cite{pierucci2014smoothing}, $f$ is non-smooth and $\cB$ is also the identity matrix. 
The proposed \spmp can be seen as a blend of the successful components resp. of the Composite Conditional Gradient algorithm~\cite{harchaoui2013conditional} and the Composite Mirror-Prox~\cite{CoMP13}, that enjoys the optimal complexity bound $O(1/\epsilon^2)$ on the total number of LMO calls, yet solves a broader class of convex problems than previously considered. 

\paragraph{Outline}
The paper is organized as follows. In Section 2, we describe the norm-regularized nonsmooth problem of interest and illustrate it with several examples. In Section 3,  we present the conditional gradient type method based on an inexact Mirror-Prox framework for structured variational inequalities. In Section 4, we present promising experimental results showing the interest of the approach in comparison to competing methods, resp. on a collaborative filtering for movie recommendation and link prediction for social network analysis applications. 
\section{Framework and assumptions}

We present here our theoretical framework, which hinges upon a smooth convex-concave saddle point reformulation of 
the norm-regularized non-smooth minimization (\ref{nonsmoothproblem}). We shall use the following notations throughout the paper.
For a given norm  $\|\cdot\|$, we define the dual norm as $\|s\|_*=\max_{\|x\|\leq 1}\langle  s, x\rangle$. For any $x\in\bR^{m\times n}$, $\|x\|_2=\|x\|_F =(\sum_{i=1}^m\sum_{j=1}^n |x_{ij}|^2)^{1/2}$. 

\paragraph{Problem} 
We consider the composite minimization problem
\begin{equation}\label{nonsmoothproblem}
\Opt =\min _{x\in X} \: f(x)+ \|\cB x\|
\end{equation}
where $X$ is a closed convex set in the  Euclidean space $E_x$; $x\mapsto \cB x$ is a linear mapping from $X$ to $Y(\supset \cB X)$, where $Y$ is a closed  convex set in the Euclidean space  $E_y$. We make two important assumptions on the function $f$ and the norm $\|\cdot\|$ defining the regularization penalty, explained below.

\paragraph{Fenchel-type Representation}
The non-smoothness of $f$ can be challenging to tackle. However, in many cases of interest, 
the function $f$ enjoys a favorable structure that allows to tackle it with smoothing techniques.  
We assume that the norm $f(x)$ is a non-smooth convex function given by 
\begin{equation}\label{eq:fenchel-rep}
f(x)=\max_{z\in Z} \: \Phi(x,z)
\end{equation}a
where $\Phi(x,z)$ is a smooth convex-concave function and $Z$ is a convex and compact set in the Euclidean space $E_z$.
Such representation was introduced and developed in~\cite{cox2013dual,CoMP13,juditsky13variational}, for the purpose 
of non-smooth optimisation. Fenchel-type representability can be interpreted as a general form of the smoothing-favorable structure 
of non-smooth functions used in the Nesterov smoothing technique~\cite{nesterov2005smooth}. Representations of this type 
are readily available for a wide family of ``well-structured'' nonsmooth functions $f$; see Sec.~\ref{sec:exps} for examples.

\paragraph{Composite Linear Minimization Oracle}
Proximal-gradient-type algorithms require the computation of a proximal operator at each iteration, i.e. 
\begin{equation}
\min_{y\in Y} \: \left\{\frac{1}{2}\|y\|_2^2+\langle \eta, y\rangle +\alpha\|y\| \right\}.
\end{equation}
For several cases of interest, described below, the computation of the proximal operator can be expensive or intractable. A classical example 
is the nuclear norm, whose proximal operator boils down to singular value thresholding, therefore requiring a full singular value decomposition. 
In contrast to the proximal operator, the linear minimization oracle can much cheaper. The linear minimization oracle (LMO) is a routine
which, given an input $\alpha>0$ and $\eta\in E_y$, returns a point 
\begin{equation} 
\min_{y\in Y} \: \left\{\langle \eta, y\rangle +\alpha\|y\| \right\}
\end{equation}
In the case of the nuclear-norm, the LMO only requires the computation of the top pair 
of eigenvectors/eigenvalues, which is an order of magnitude fast in time-complexity. 

\paragraph{Saddle Point Reformulation.} 
The crux of our approach is a smooth convex-concave saddle point reformulation of (\ref{nonsmoothproblem}). After massaging the 
saddle-point reformulation, we consider the variational inequality associated with the obtained saddle-point problem. 
For a constrained smooth optimisation problem, the corresponding variational inequality provides
the sufficient and necessary condition for an optimal solution to the problem~\cite{Bauschke:2011,bertsekas:2015}. For non-smooth optimization problems,
the corresponding variational inequality is directly related to the accuracy certificate used to guarantee
the accuracy of a solution to the optimisation problem; see Sec. 2.1 in~\cite{CoMP13} and~\cite{Nem10}. We shall present then an algorithm to solve the variational inequality established below, that leverages its particular structure.

Assuming that $f$ admits a Fenchel-type representation  (\ref{eq:fenchel-rep}), we rewrite (\ref{nonsmoothproblem}) in epigraph form
\begin{equation*}\label{nonsmoothrewrite}
\min _{x\in X, y\in Y,\tau\geq \|y\|}\max_{z\in Z} \: \left\{\Phi(x,z)+ \tau: y=\cB x\right\},
\end{equation*}
which, with a properly selected $\rho>0$, can be further approximated by 
\begin{align}\label{saddle problem}
\hat\Opt&=\min_{x\in X, y\in Y, \tau\geq\|y\|}\max_{z\in Z} \: \left \{\Phi(x,z)+\tau +\rho\|y-\cB x\|_2\right\}\\
&=\min_{x\in X, y\in Y, \tau\geq\|y\|}\max_{z\in Z, \|w\|_2\leq 1} \: \left \{\Phi(x,z)+\tau +\rho\langle y-\cB x, w\rangle\right\} \; .
\end{align}
In fact, when $\rho$ is large enough
one can always guarantee $\hat\Opt=\Opt$. It is indeed sufficient to set $\rho$ as the Lipschitz constant of $\|\cdot\|$ with respect to $\|\cdot\|_2$.

Introduce the variables $u:=[x,y;z,w]$ and $v:=\tau$. The variational inequality associated with the above saddle point problem is fully described by the domain
\begin{equation*}
\begin{array}{rcl}
X_+&=&\{x_+=[u;v]: x\in X, y\in Y, z\in Z, \|w\|_2\leq 1,\tau\geq \|y\|\}
\end{array}\end{equation*}
and the monotone vector field
\begin{equation*}
F(x_+=[u;v])=[F_u(u);F_v] \; ,
\end{equation*}
where 
\begin{align*}
\label{semi-struct-var-ineq}
F_u\left(u=\begin{bmatrix} x\\y\\z\\w\end{bmatrix}\right)&=
\begin{bmatrix}  
\nabla_x \Phi(x,z)-\rho\cB^Tw\\
\rho w\\
-\nabla_z\Phi(x,z)\\
\rho(\cB x-y)
 \end{bmatrix}, &\quad  F_v(v=\tau)&=1.
\end{align*}
In the next section, we present an efficient algorithm to solve this type of variational inequality,
which enjoys a particular structure; we call such an inequality \emph{semi-structured}.
\section{\Spmp for Semi-structured Variational Inequalities} 

Semi-structured variational inequalities (Semi-VI) enjoy a particular product structure, that allows to get the best 
of two worlds, namely the proximal setup (where the proximal operator can be computed)
and the LMO setup (where the linear minimization oracle can be computed). Basically, 
the domain $X$ is decomposed as a Cartesian product over two sets $X = X_1 \times X_2$,
such that $X_1$ admits a proximal-mapping while $X_2$ admits a linear minimization oracle. 
We now describe the main theoretical and algorithmic components of the \spmp algorithm, 
resp. in Sec.~\ref{sec:CMP} and in Sec.~\ref{sec:CCG}, and finally describe the overall 
algorithm in Sec.~\ref{sec:semivi-algo}. 

\subsection{Composite Mirror-Prox with Inexact Prox-mappings}\label{sec:CMP}

We first present a new algorithm, which can be seen as an extension of the composite Mirror Prox algorithm, denoted CMP for brevity, 
that allows inexact computation of the Prox-mappings, and can solve a broad class of variational inequalites.
The original Mirror Prox algorithm was introduced in~\cite{Nem04}, and was extended to composite minimization
in~\cite{CoMP13} assuming exact computations of Prox-mappings.

\paragraph{Structured Variational Inequalities.} We consider the variational inequality VI$(X,F)$:
$$ \text{Find } x_*\in X: \langle F(x),x-x_*\rangle\geq 0, \forall x\in X$$
with domain $X$ and operator $F$ that satisfy the assumptions ({\bf A}.1)--({\bf A}.4) below.
\begin{enumerate}[({\textbf{A}}.1)]
\item Set $X\subset E_u\times E_v$ is closed convex and its projection $PX=\{u: x=[u;v]\in X\}\subset U$, where $U$ is convex and closed, $E_u,E_v$ are Euclidean spaces;
\item The function $\omega(\cdot):U\to \bR$ is continuously differentiable 
and also 1-strongly convex w.r.t. some norm\footnote{There is a slight abuse of notation here. The norm here is not the same as the one in problem (\ref{nonsmoothproblem})}
 $\|\cdot\|$.
This defines the Bregman distance
\begin{equation*}
V_u(u')=\omega(u')-\omega(u)-\langle\omega'(u),u'-u\rangle\geq\frac{1}{2}\|u'-u\|^2 \; .
\end{equation*}

\item The operator $F(x=[u,v]): X\to E_u\times E_v$ is monotone and of form $ F(u,v)=[F_u(u);F_v]$
with $F_v\in E_v$ being a constant and $F_u(u)\in E_u$ satisfying the condition
\begin{equation*}
\begin{array}{c}
\forall u,u'\in U: \|F_u(u)-F_u(u')\|_*\leq L\|u-u'\|+M
\end{array}\end{equation*}
for some $L<\infty, M<\infty$;

\item The linear form $\langle F_v,v\rangle$ of $[u;v]\in E_u\times E_v$ is bounded from below on $X$ and is coercive on $X$ w.r.t. $v$: whenever $[u^t;v^t]\in X$, $t=1,2,...$ is a sequence such that $\{u^t\}_{t=1}^\infty$ is bounded and $\|v^t\|_2\to\infty$ as $t\to\infty$, we have $\langle F_v,v^t\rangle \to\infty$, $t\to\infty$.
\end{enumerate}


\paragraph{$\epsilon$-Prox-mapping}
In the Composite Mirror Prox with exact Prox-mappings~\cite{CoMP13}, the quality of an iterate, in the course
of the algorithm, is measured through the so-called dual gap function 
\begin{equation*}
\epsilonvi(x\big|X,F) = \sup_{y\in X} \; \langle F(y),x-y\rangle \; .
\end{equation*}
We give in Appendix~\ref{sect:preliminaries} a refresher on dual gap functions, for the reader's convenience. 
We shall establish the complexity bounds in terms this dual gap function for our algorithm, which directly provides an accuracy certificate
along the iterations. However, we first need to define what we mean by an inexact prox-mapping. Inexact proximal mapping
were recently considered in the context of accelerated proximal gradient algorithms~\cite{inexact-prox-grad}. The definition we give below
is more general, allowing for non-Euclidean proximal-mappings. 

We introduce here the notion of $\epsilon$-prox-mapping ($\epsilon\geq 0)$. 
For $\xi=[\eta;\zeta]\in E_u\times E_v$ and $x=[u;v]\in X$, let us define the subset $P_x^\epsilon(\xi)$ of $X$ as
\begin{equation*}\label{eq:epsprox}
P_x^\epsilon(\xi)=\{\hat x=[\hat u;\hat v]\in X: \langle \eta+\omega'(\hat u)-\omega'(u),\hat u-s\rangle+\langle \zeta,\hat v-w\rangle \leq\epsilon\,\,\forall [s;w]\in X\}.
\end{equation*}
When $\epsilon=0$, this reduces to the exact prox-mapping, in the usual setting,
that is 
\begin{equation*}
P_x(\xi)=\Argmin_{[s;w]\in X} \: \left\{\langle \eta,s\rangle+\langle\zeta,w\rangle +V_u(s)\right\} \; .
\end{equation*}
When $\epsilon>0$, this yields our definition of an inexact prox-mapping, with inexactness parameter $\epsilon$. 
Note that for any $\epsilon\geq 0$, the set $P_x^\epsilon(\xi=[\eta;\gamma F_v])$ is well defined whenever $\gamma>0$. The Composite Mirror-Prox with Inexact Prox-mappings is outlined in Algorithm~\ref{CMP}.

\begin{algorithm}[H]
\caption{Composite Mirror Prox Algorithm (CMP) for VI$(X,F)$}
\label{CMP}
\begin{algorithmic}
\STATE \textbf{Input:} stepsizes $\gamma_t>0$, inexactness $\epsilon_t\geq 0$, $t = 1,2, \ldots$
\STATE Initialize $x^1=[u^1;v^1]\in X$
\FOR{$t=1,2,\ldots, T$}
\STATE 
\setlength\belowdisplayskip{2pt}
\begin{equation}\label{MirrorProx}
\begin{array}{rcl}
y^{t}:=[\hat{u}^{t};\hat{v}^{t}]&\in&P_{x^t}^{\epsilon_t}(\gamma_t F(x^t))=P_{x^t}^{\epsilon_t}(\gamma_t[ F_u(u^t);F_v])\\
x^{t+1}:=[u^{t+1};v^{t+1}]&\in&P_{x^t}^{\epsilon_t}(\gamma_t F(y^t))=P_{x^t}^{\epsilon_t}(\gamma_t[ F_u(\hat{u}^t);F_v])
\end{array}
\end{equation}
\ENDFOR\\
\STATE \textbf{Output:} $\overline x_T:=[\bar{u}_T;\bar{v}_T] ={(\sum_{t=1}^T\gamma_t )}^{-1}{\sum_{t=1}^T \gamma_t y^t}$
\end{algorithmic}
\end{algorithm}

Note that this composite version of Mirror Prox algorithm works essentially {\sl as if} there were no $v$-component at all. Therefore, the proposed 
algorithm is a not-trivial extension of the Composite Mirror-Prox with \emph{exact prox-mappings}, both from a theoretical and algorithmic
point of views. We establish below the theoretical convergence rate; see Appendix for the proof. 

\begin{theorem}\label{thm:theMP}
Assume that the sequence of step-sizes $(\gamma_t)$ in the CMP algorithm satisfy
\begin{equation}\label{gammaupperbound}
\sigma_t:=\gamma_t\langle F_u(\hat{u}^t)-F_u(u^t),\hat{u}^t-u^{t+1}\rangle-V_{\hat{u}^t}(u^{t+1})-V_{u^t}(\hat{u}^{t})\le \gamma_t^2M^2 \, , \quad t=1,2,\ldots, T \; .
\end{equation}
Then, denoting $\Theta[X]=\sup_{[u;v]\in X}V_{u^1}(u)$, for a sequence of inexact prox-mappings with inexactness $\epsilon_t\geq 0$, we have
\begin{equation}\label{epsilonvismall}
\epsilonvi(\bar{x}_T\big|X,F):=\sup_{x\in X} \; \left\langle F(x), \bar{x}_T-x\right\rangle \leq \frac{\Theta[X]+M^2{\sum}_{t=1}^T\gamma_t^2+2{\sum}_{t=1}^T\epsilon_t}{\sum_{t=1}^T\gamma_t}.
\end{equation}
\end{theorem}

\textbf{Remarks }
Note that the assumption on the sequence of step-sizes $(\gamma_t)$ is clearly satisfied when $\gamma_t\leq ({\sqrt{2}L})^{-1}$. When $M=0$, it is satisfied as long as $\gamma_t\leq L^{-1}$. 

\begin{corollary}\label{cor:theMP}
Assume further that  $X=X_1\times X_2$, and let $F$ be the monotone vector field associated with the saddle point problem 
\begin{equation}\label{saddlepoint}
\SadVal=\min_{x^1\in X_1}\max_{x^2\in X_2}\Phi(x^1,x^2),
\end{equation} two induced convex optimization problems
\begin{equation}\label{primaldual}
\begin{array}{rclr}
\Opt(P)&=&\min_{x^1\in X_1}\left[\overline{\Phi}(x^1)=\sup_{x^2\in X_2}\Phi(x^1,x^2)\right]&\qquad\qquad(P)\\
\Opt(D)&=&\max_{x^2\in X_2}\left[\underline{\Phi}(x^2)=\inf_{x^1\in X_1}\Phi(x^1,x^2)\right]&\qquad\qquad(D)\\
\end{array}
\end{equation}
with convex-concave locally Lipschitz continuous cost function $\Phi$. In addition, assuming  that problem $(P)$ in {\rm (\ref{primaldual})} is solvable with optimal solution $x^1_*$ and denoting by $\bar{x}_T^1$ the projection of $\bar{x}_T\in X=X_1\times X_2$ onto $X_1$, we have
\be
\overline{\Phi}(\bar{x}_T^1)-\Opt(P)\leq \left[{\sum}_{t=1}^T\gamma_t\right]^{-1}\left[\Theta[\{x^1_*\}\times X_2]+M^2{\sum}_{t=1}^T\gamma_t^2+2{\sum}_{t=1}^T\epsilon_t\right].
\ee{optimalitygap}
\end{corollary}

The theoretical convergence rate established in Theorem~\ref{thm:theMP} and Corollary~\ref{cor:theMP} generalizes the previous result established in Corollary 3.1 in~\cite{CoMP13}
for CMP with exact prox-mappings. Indeed, when exact prox-mappings are used, we recover the result of~\cite{CoMP13}. When inexact prox-mappings are used, the errors due to the inexactness of the prox-mappings accumulates and is reflected in the bound (\ref{epsilonvismall}) and (\ref{optimalitygap}). 


\subsection{Composite Conditional Gradient}\label{sec:CCG}
We now turn to a variant of the composite conditional gradient algorithm, denoted CCG, tailored for a particular class of problems,
which we call \emph{smooth semi-linear problems}. The composite conditional gradient algorithm 
was introduced in~\cite{harchaoui2013conditional}. We present an extension here which will turn 
to be especially tailored for sub-problems that will be solved in Sec.~\ref{sec:semivi-algo}. 

\paragraph{Minimizing Smooth Semi-linear Problems.} We consider the smooth semi-linear problem
\begin{equation}\label{CCGProblem}
\min_{x=[u;v]\in X} \left\{\phi^+(u,v)=\phi(u)+\langle\theta,v\rangle\right\}
\end{equation}
represented by the pair $(X;\phi^+)$ such that the following assumptions are satisfied. We assume that
\begin{enumerate}[i)]
\item 
$X\subset E_u\times E_v$ is closed convex and its projection $PX\subset U$, where $U$ is convex and compact;
\item 
$\phi(u):U\to\bR$ be a convex continuously differentiable function, and there exists $1< \kappa \leq 2$ and $L<\infty$ such that
\begin{equation}\label{CCG1}
\phi(u')\leq\phi(u)+\langle \nabla\phi(u),u'-u\rangle + {L_0\over\kappa}\|u'-u\|^{\kappa}\,\,\forall u,u'\in U ;
\end{equation}

\item $\theta\in E_v$ be such that every linear function on $E_u\times E_v$ of the form
\begin{equation}
[u;v]\mapsto \langle \eta,u\rangle + \langle \theta,v\rangle 
\end{equation}
with $\eta\in E_u$ attains its minimum on $X$ at some point $x[\eta]=[u[\eta];v[\eta]]$; we have at our disposal a \emph{Composite Linear Minimization Oracle} (LMO) which, given on input $\eta\in E_u$, returns $x[\eta]$.
\end{enumerate}

\begin{algorithm}
\caption{Composite Conditional Gradient Algorithm {\bf CCG$(X,\phi(\cdot),\theta;\epsilon)$}}
\label{CCG}
\begin{algorithmic}
\STATE \textbf{Input:} accuracy $\epsilon>0$ and $\gamma_t = 2/(t+1), t = 1,2,\ldots$
\STATE Initialize $x^1=[u^1;v^1]\in X$ and 
\FOR{$t=1,2,\ldots $}
\STATE Compute $\delta_t=\langle g_t,u^t-u^t[g_t]\rangle+\langle \theta,v^t-v^t[g_t]\rangle$, where $g_t=\nabla\phi(u^t)$;
\IF {$\delta_t\leq \epsilon$}
\STATE Return $x^t=[u^t;v^t]$ 
\ELSE
\STATE Update $x^{t+1}=[u^{t+1};v^{t+1}]\in X$ such that
$
\phi^+(x^{t+1})\leq \phi^+\left(x^t+\gamma_t(x^t[g_t]-x^t)\right)
$
\ENDIF
\ENDFOR\\
\end{algorithmic}
\end{algorithm}


The algorithm is outlined in Algorithm~\ref{CCG}. Note that CCG works essentially {\sl as if} there were no $v$-component at all.
The CCG algorithm enjoys a convergence rate in $O(t^{-(\kappa-1)})$ in the evaluations of the function $\phi^+$, 
and the accuracy certificates $(\delta_t)$ enjoy the same rate $O(t^{-(\kappa-1)})$ as well, for solving problems of type~(\ref{CCGProblem}). 
See Appendix for details and the proof.

\begin{proposition}\label{prop:theCCG}
Denote $D$ the $\|\cdot\|$-diameter of $U$. When solving problems of type~(\ref{CCGProblem}), the sequence of iterates $(x^t)$ of CCG satisfies
\begin{equation}\label{suboptimality}
\epsilon_t:=\phi^+(x^t)-\displaystyle\min_{x\in X}\phi^+(x)\leq {2L_0D^\kappa\over\kappa(3-\kappa)}\left(\frac{2}{t+1}\right)^{\kappa-1},\,t\geq2
\end{equation}
In addition, the accuracy certificates $(\delta_t)$ satisfy
\begin{equation}\label{certificate}
\min_{1\leq s\leq t} \; \delta_s \leq O(1)L_0D^\kappa \left(\frac{2}{t+1}\right)^{\kappa-1},\,t\geq2.
\end{equation}
\end{proposition}

\subsection{\Spmp for Semi-structured Variational Inequality}
\label{sec:semivi-algo}
We now give the full description of a special class of variational inequalities, called \emph{semi-structured variational inequalities}. 
This family of problems encompasses both cases that we discussed so far in Section \ref{sec:CMP} and \ref{sec:CCG}. But most importantly, it also covers many other problems that do not fall into these two regimes and in particular, our essential problem of interest (\ref{nonsmoothproblem}).

\paragraph{Semi-structured Variational Inequalities.} 
The class of semi-structured variational inequalities allows to go beyond Assumptions $(\textbf{A}.1)-(\textbf{A}.4)$, by assuming more structure. This structure is consistent with what we call a \emph{semi-proximal} setup, which encompasses both the regular \emph{proximal setup} and the regular
\emph{linear minimization setup} as special cases. Indeed, we consider a class of
variational inequality VI$(X,F)$ that satisfies, in addition to Assumptions $(\textbf{A}.1)-(\textbf{A}.4)$, the following assumptions:
\begin{enumerate}[({\textbf{S}}.1)]
\item \emph{Proximal setup for $X$}:  we assume that $E_u=E_{u_1}\times E_{u_2}$, $E_v=E_{v_1}\times E_{v_2}$, and  $U\subset U_1\times U_2$, $X=X_1\times X_2$ with $X_i\in E_{u_i}\times E_{v_i}$ and $P_iX=\{u_i:[u_i;v_i]\in X_i\}\subset U_i$ for $i=1,2$, where $U_1$ is convex and closed,  $U_2$ is convex and compact. We also assume that $\omega(u)=\omega_1(u_1)+\omega_2(u_2)$ and $\|u\|=\|u_1\|_{E_{u_1}}+\|u_2\|_{E_{u_2}}$, with $\omega_2(\cdot): U_2\to\bR$ continuously differentiable such that  
\begin{equation*}
\omega_2(u_2')\leq\omega_2(u_2)+\langle\nabla \omega_2(u_2), u_2'-u_2\rangle+ \frac{L_0}{\kappa} \|u_2'-u_2\|_{E_{u_2}}^\kappa,\forall u_2,u_2'\in U_2;
\end{equation*}
for a particular $1< \kappa \leq 2$ and $L_0<\infty$. Furthermore, we assume that the $\|\cdot\|_{E_{u_2}}$-diameter of $U_2$ is bounded by some $D>0$..

\item \emph{Proximal mapping on $X_1$}: we assume that for any $\eta_1\in E_{u_1}$ and $\alpha>0$,  we have at disposal easy-to-compute prox-mappings of the form,
\begin{equation*}
\Prox_{\omega_1}(\eta_1,\alpha):=\min_{x_1=[u_1;v_1]\in X_1}\; \left\{\omega_1(u_1)+ \langle\eta_1, u_1\rangle+\alpha\langle F_{v_1}, v_1\rangle\right\}.
\end{equation*}

\item \emph{Linear minimization on $X_2$}: we assume that we we have at our disposal Composite Linear Minimization Oracle (LMO), which given any input $\eta_2\in E_{u_2}$ and $\alpha>0$, returns an optimal solution to the minimization problem with linear form, that is, 
\begin{equation*}
\LMO(\eta_2,\alpha):=\min_{x_2=[u_2;v_2]\in X_2}\;\left\{ \langle \eta_2, u_2\rangle+\alpha\langle F_{v_2}, v_2\rangle\right\}.
\end{equation*}
\end{enumerate}

\paragraph{Semi-proximal setup}
We denote such problems as Semi-VI$(X,F)$. On the one hand, when $U_2$ is a singleton, we get the \emph{full-proximal setup}. On the other hand, when $U_1$ is a singleton, we get the \emph{full linear-minimization-oracle setup} (full LMO setup). In the gray zone in between, we get the \emph{semi-proximal setup}.

\paragraph{The \spmp algorithm.}  We finally present here our main contribution, the \spmp algorithm, 
which solves the semi-structured variational inequality under $(\textbf{A}.1)-(\textbf{A}.4)$ and $(\textbf{S}.1)-(\textbf{S}.3)$. The \spmp algorithm blends both CMP and CCG. Basically, for sub-domain $X_2$ given by LMO, instead of computing exactly the prox-mapping, 
we mimick inexactly the prox-mapping via a conditional gradient algorithm in the composite Mirror Prox algorithm. 
For the sub-domain $X_1$, we compute the prox-mapping as it is.

\paragraph{Course of the \spmp algorithm}
Basically, at step $t$, we first update $y_1^t=[\hat u_1^t;\hat v_1^t]$ by computing the exact prox-mapping and update
$y_2^t=[\hat u_2^t; \hat v_2^t]$  by running the composite conditional gradient algorithm to problem 
(\ref{CCGProblem}) specifically with 
$$X=X_2, \phi(\cdot)=\omega_2(\cdot)+\langle \gamma_tF_{u_2}(u_2^t)-\omega_2'(u_2^t),\cdot \rangle , \text{ and } \theta=\gamma_t F_{v_2},$$ 
until $
\delta(y_2^t)=\max_{y_2\in X_2} \langle \nabla\phi^+(y_2^t), y_2^t-y_2\rangle \leq \epsilon_t$. 
We then update $x_1^{t+1}=[u_1^{t+1};v_1^{t+1}]$  and $x_2^{t+1}=[u_2^{t+1};v_2^{t+1}]$ similarly except this time taking the value of the operator at point $y^t$. Combining the results in Theorem \ref{thm:theMP} and Proposition \ref{prop:theCCG}, we arrive at the following complexity bound.

\begin{algorithm}[H]
\caption{\textbf{\spmpnospace}  Algorithm for Semi-VI$(X,F)$}
\label{CMP_CCG}
\begin{algorithmic}
\STATE \textbf{Input:} stepsizes $\gamma_t>0$, accuracies $\epsilon_t\geq 0$, $t = 1,2, \ldots$
\STATE[1] Initialize $x^1=[x_1^1;x_2^1]\in X$, where $x_1^1=[u_1^1;v_1^1];x_2^1=[u_2^1,;v_2^1]$.
\FOR{$t=1,2,\ldots, T$}
\STATE[2] Compute $y^t=[y_1^t;y_2^t]$ that
\begin{equation*}\label{MirrorProx1}
\begin{array}{rcl}
y_1^{t}:=[\hat{u}_1^{t};\hat{v}_1^{t}]&=&\Prox_{\omega_1}(\gamma_t F_{u_1}(u_1^t)-\omega_1'(u_1^t),\gamma_t)\\
y_2^{t}:=[\hat{u}_2^{t};\hat{v}_2^{t}]&=&\textbf{CCG}(X_2, \omega_2(\cdot)+\langle \gamma_tF_{u_2}(u_2^t)-\omega_2'(u_2^t),\cdot \rangle, \gamma_tF_{v_2};\epsilon_t)\\
\end{array}
\end{equation*}

\STATE[3] Compute $x^{t+1}=[x_1^{t+1};x_2^{t+1}]$ that
\begin{equation*}\label{MirrorProx2}
\begin{array}{rcl}
x_1^{t+1}:=[u_1^{t+1};v_1^{t+1}]&=&\Prox_{\omega_1}(\gamma_t F_{u_1}(\hat{u}_1^t)-\omega_1'(u_1^t),\gamma_t)\\
x_2^{t+1}:=[u_2^{t+1};v_2^{t+1}]&=&\textbf{CCG}(X_2, \omega_2(\cdot)+\langle \gamma_tF_{u_2}(\hat{u}_2^t)-\omega_2'(u_2^t),\cdot \rangle, \gamma_tF_{v_2};\epsilon_t)\\
\end{array}
\end{equation*}
\ENDFOR\\
\STATE \textbf{Output:} $\overline x_T:=[\bar{u}_T;\bar{v}_T] ={(\sum_{t=1}^T\gamma_t )}^{-1}{\sum_{t=1}^T \gamma_t y^t}$
\end{algorithmic}
\end{algorithm}

\begin{proposition}\label{prop:thecomplexity}
Under the assumption $(\textbf{A}.1)-(\textbf{A}.4)$ and $(\textbf{S}.1)-(\textbf{S}.3)$ with $M=0$,  for the  outlined algorithm to return an $\epsilon$-solution to the variational inequality $VI(X,F)$, the total number of Mirror Prox steps required does not exceed
\begin{equation*}
\text{Total number of steps} = O\left(\frac{L\Theta[X]}{\epsilon}\right)
\end{equation*}
and the total number of calls to the Linear Minimization Oracle does not exceed 
\begin{equation*}
\cN =O(1)\left(\frac{L_0L^\kappa D^\kappa}{\epsilon^\kappa}\right)^{\frac{1}{\kappa-1}}\Theta[X].
\end{equation*}
In particular, if we use Euclidean proximal setup on $U_2$ with $\omega_2(\cdot)=\frac{1}{2}\|x_2\|^2$, which leads to $\kappa =2$ and $L_0=1$, then the number of LMO calls does not exceed $\cN=O(1)\left(L^2D^2(\Theta[X_1]+D^2\right)/\epsilon^2$.
\end{proposition}

\paragraph{Discussion}
The proposed \spmp algorithm enjoys the \emph{optimal complexity bounds}, 
i.e. $O(1/\epsilon^2)$, in the number of calls to LMO; see~\cite{lan2013complexity}
for the optimal complexity bounds for general non-smooth optimisation with LMO. 
Furthermore, \spmp generalizes
previously proposed approaches and improves upon them in special cases of problem (\ref{nonsmoothproblem}); see Appendix.
\section{Experiments}\label{sec:exps}
We present here illustrations of the proposed approach. We report the experimental results obtained with the proposed \spmpnospace, denoted \textbf{Semi-MP} here, and state-of-the-art competing optimization algorithms. We consider three different models, all with a non-smooth loss function and a nuclear-norm regularization penalty: i) matrix completion with $\ell_2$ data fidelity term; ii) robust collaborative filtering for movie recommendation; iii) link prediction for social network analysis. For i) \& ii), we compare to two competing approaches: a) smoothing conditional gradient proposed in \cite{pierucci2014smoothing} (denoted \text{Smooth-CG}); b) smoothing proximal gradient (\cite{Nesterov07a, Chen12}) equipped semi-proximal setup (\text{Semi-SPG}). For iii), we compare to Semi-LPADMM, using~\cite{Lan2014}, and solving proximal mapping through conditional gradient routines. Additional experiments and implementation details are given in Appendix~\ref{sec:compet-details}.

\paragraph{Matrix completion on synthetic data} We consider the matrix completion problem, with a nuclear-norm regularisation penalty and an $\ell_2$ data-fidelity term. 
We first investigate the convergence patterns of our Semi-MP and Semi-SPG under two different strategies of the inexactness, a) fixed inner CG steps and b) decaying $\epsilon_t=c/t$ as the theory suggested.
The plots in Fig.~\ref{fig:synthetic} indicate that using the second strategy with $O(1/t)$ decaying inexactness provides better and more reliable performance than using fixed number of inner steps. Similar trends are observed for the Semi-SPG. One can see that these two algorithms based on inexact proximal mappings are notably faster than applying conditional gradient on the smoothed problem. 

\begin{figure} 
  \begin{minipage}[t]{.25\textwidth}
    \includegraphics[scale=0.3]{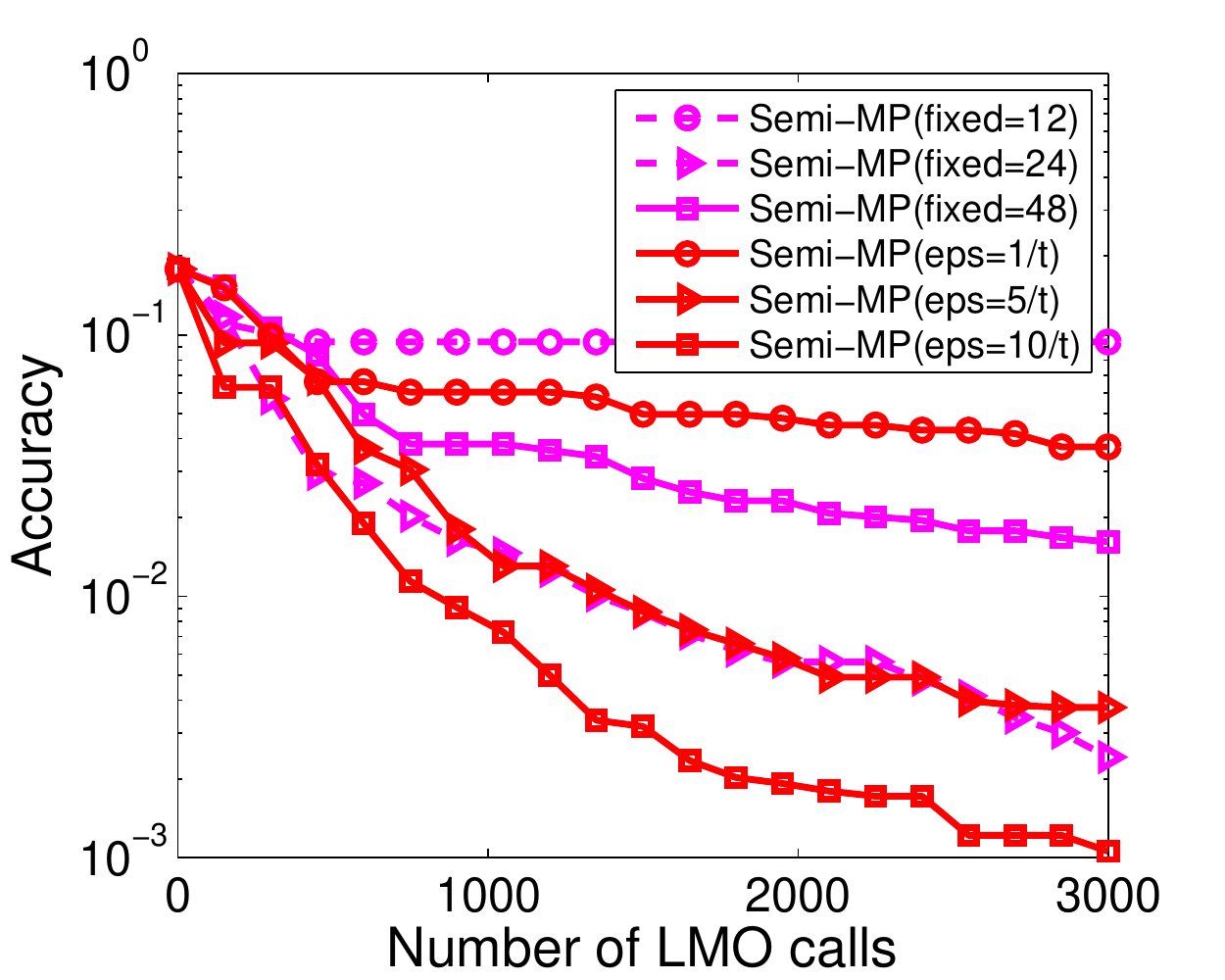}
  \end{minipage}%
  \begin{minipage}[t]{.25\textwidth}
    \includegraphics[scale=0.3]{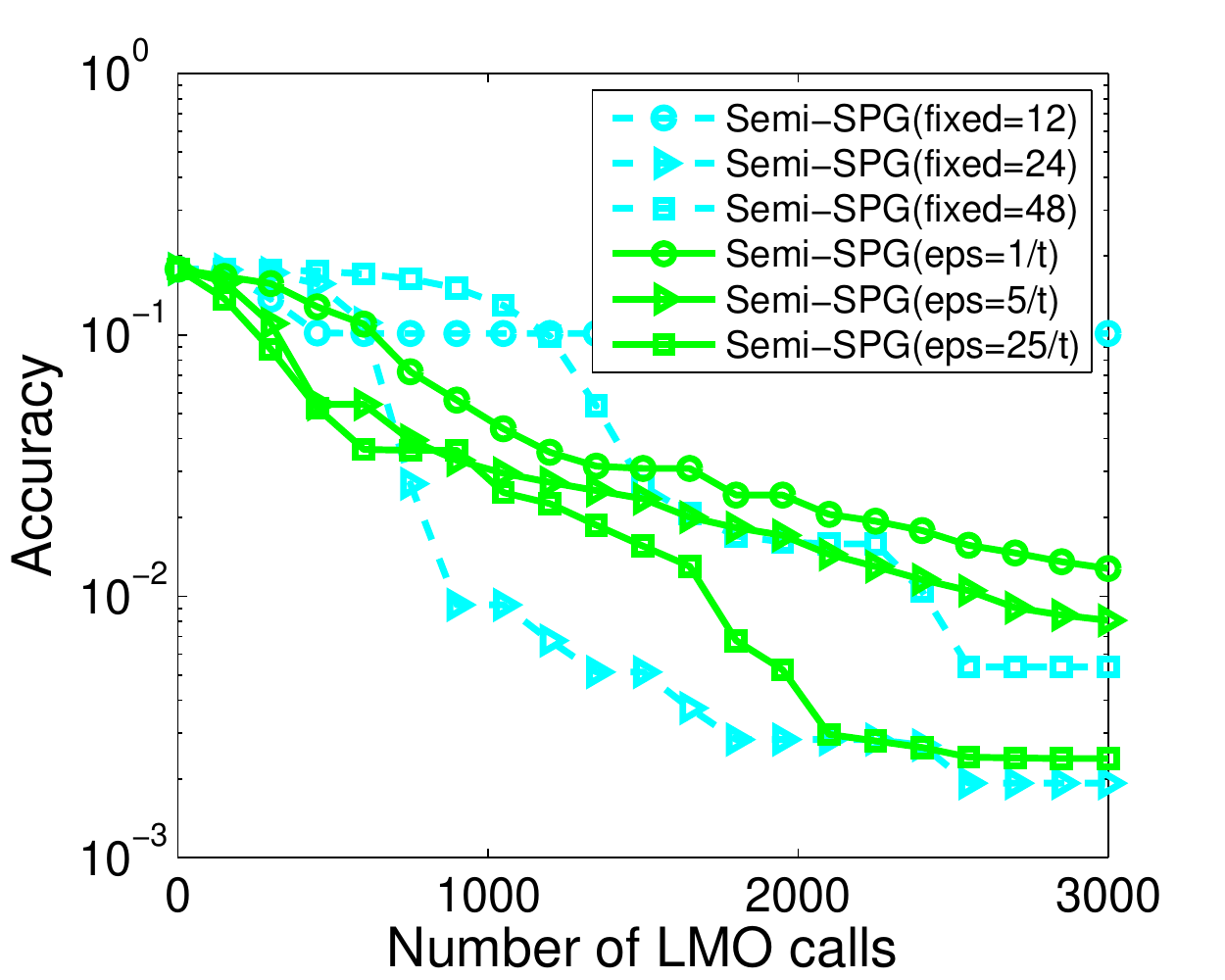}
  \end{minipage}%
  \begin{minipage}[t]{.25\textwidth}
    \includegraphics[scale=0.3]{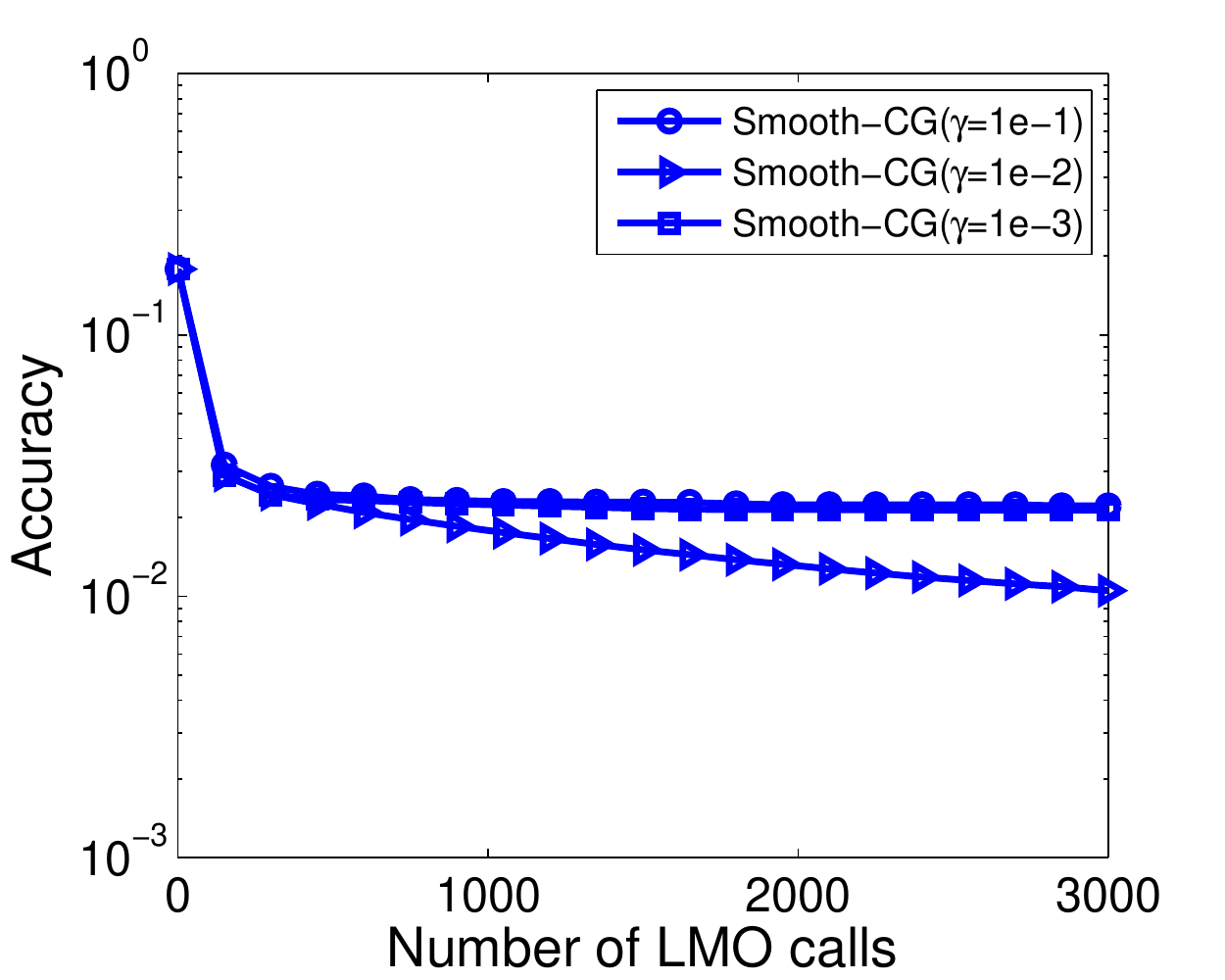}
  \end{minipage}%
   \begin{minipage}[t]{.25\textwidth}
    \includegraphics[scale=0.3]{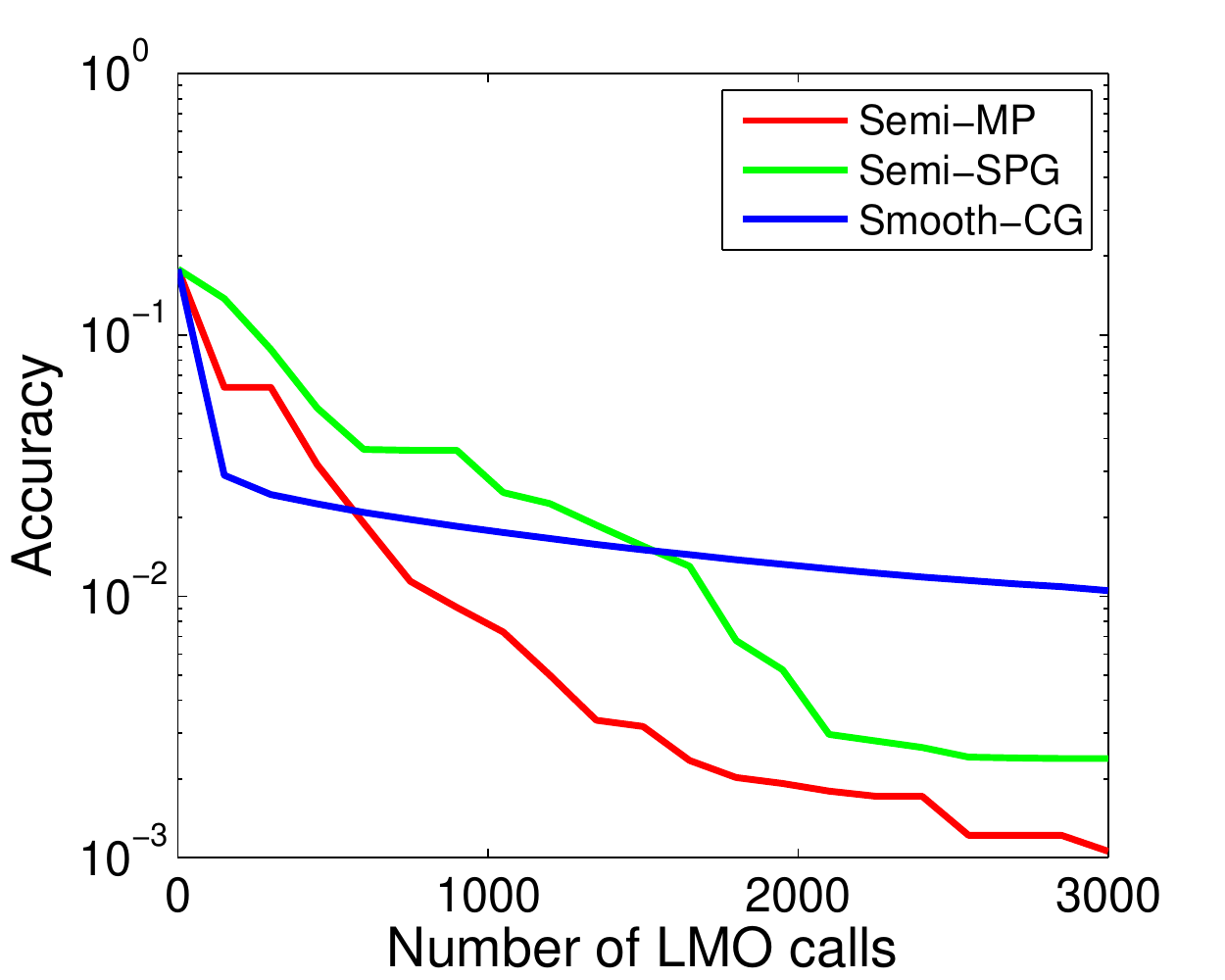}
  \end{minipage}
\caption{Matrix completion on synthetic data ($1024\times1024$): optimality gap vs the LMO calls.\\ 
   From left to right:  (a) Semi-MP; (b) Semi-SPG ; (c) Smooth-CG; (d) best of three algorithms.\\[-16pt]}
   \label{fig:synthetic}
\end{figure}

\begin{figure} 
  \begin{minipage}[t]{.25\textwidth}
    \includegraphics[scale=0.28]{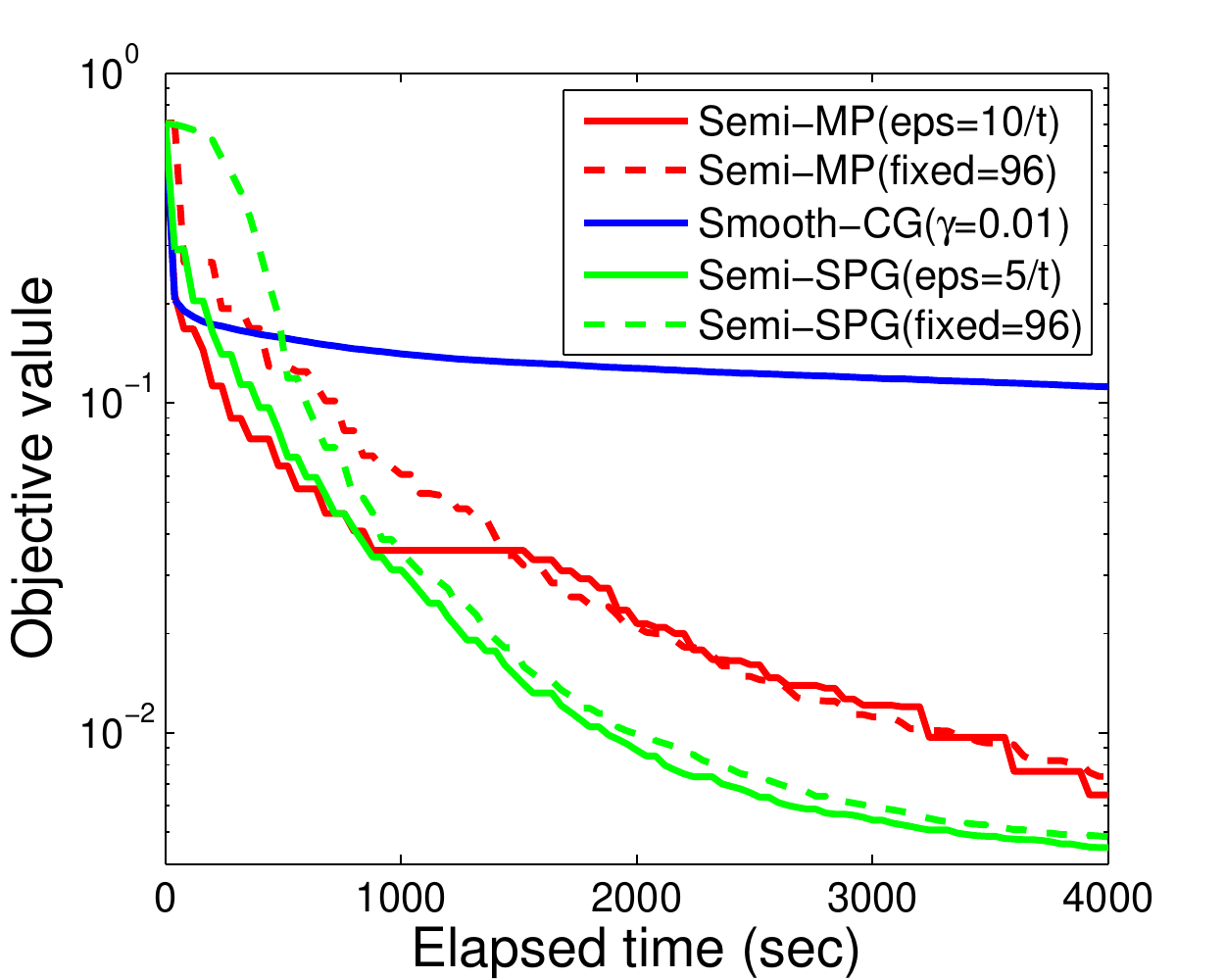}
  \end{minipage}%
  \begin{minipage}[t]{.25\textwidth}
    \includegraphics[scale=0.28]{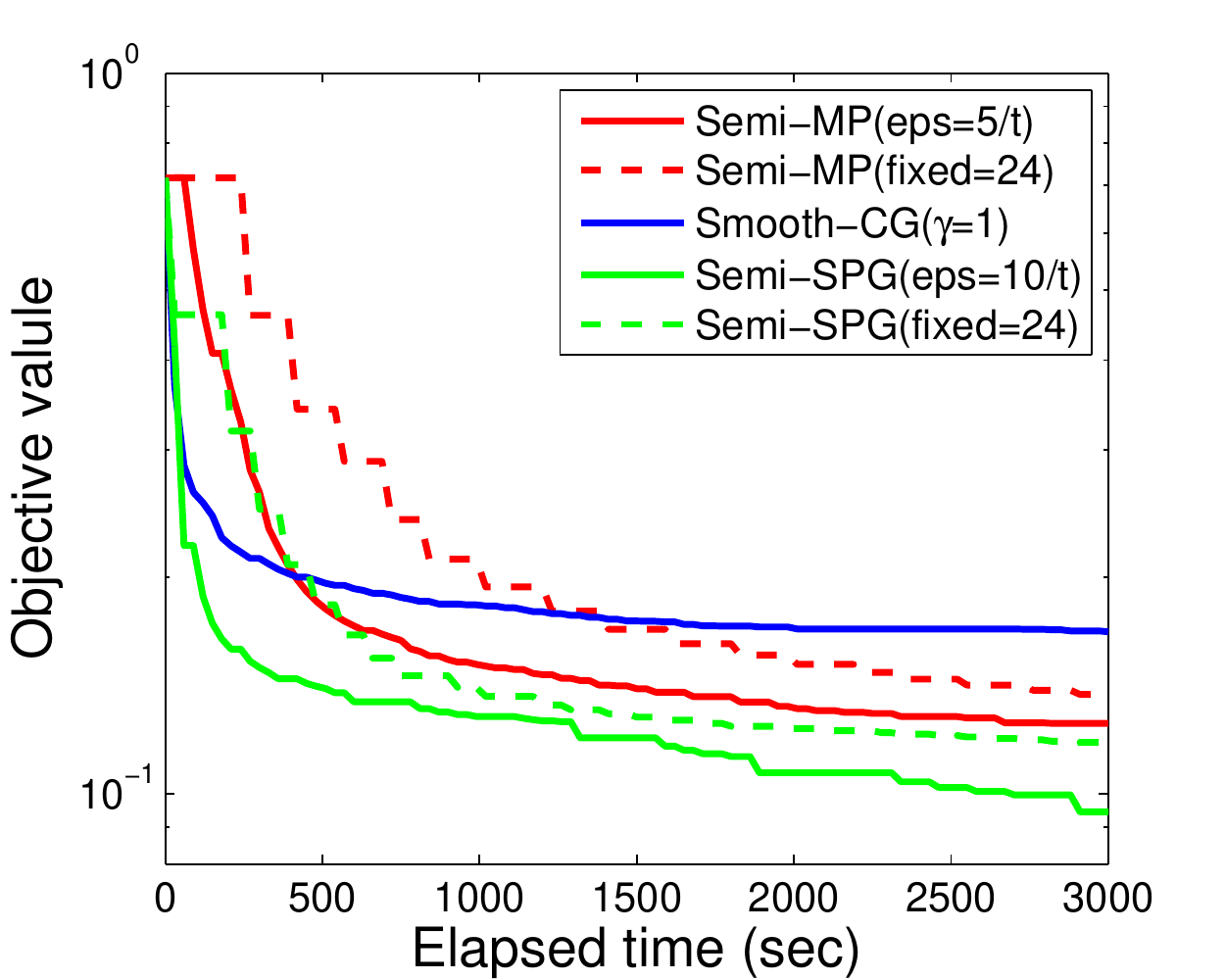}
  \end{minipage}%
  \begin{minipage}[t]{.25\textwidth}
    \includegraphics[scale=0.28]{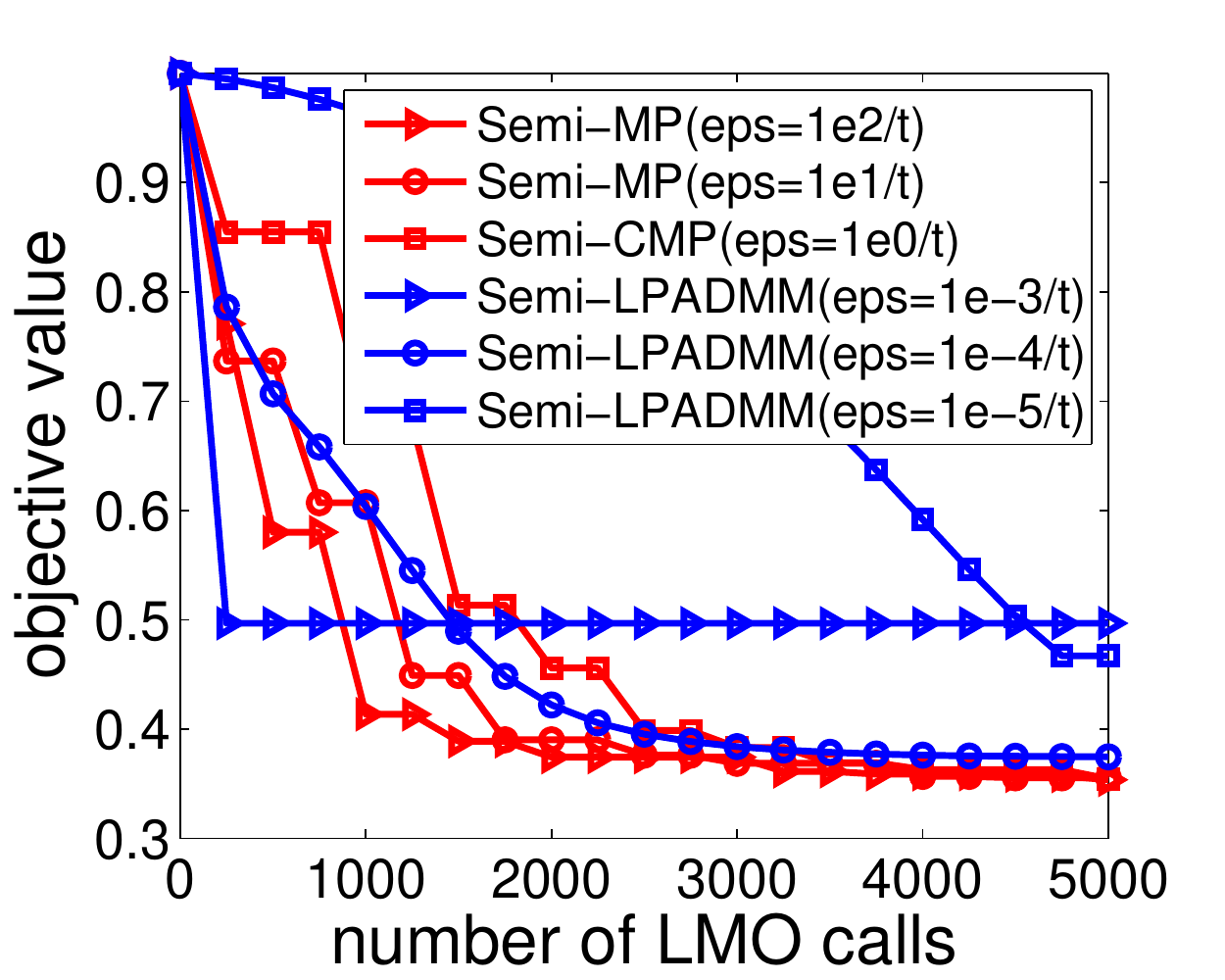}
  \end{minipage}%
 \begin{minipage}[t]{.25\textwidth}
    \includegraphics[scale=0.28]{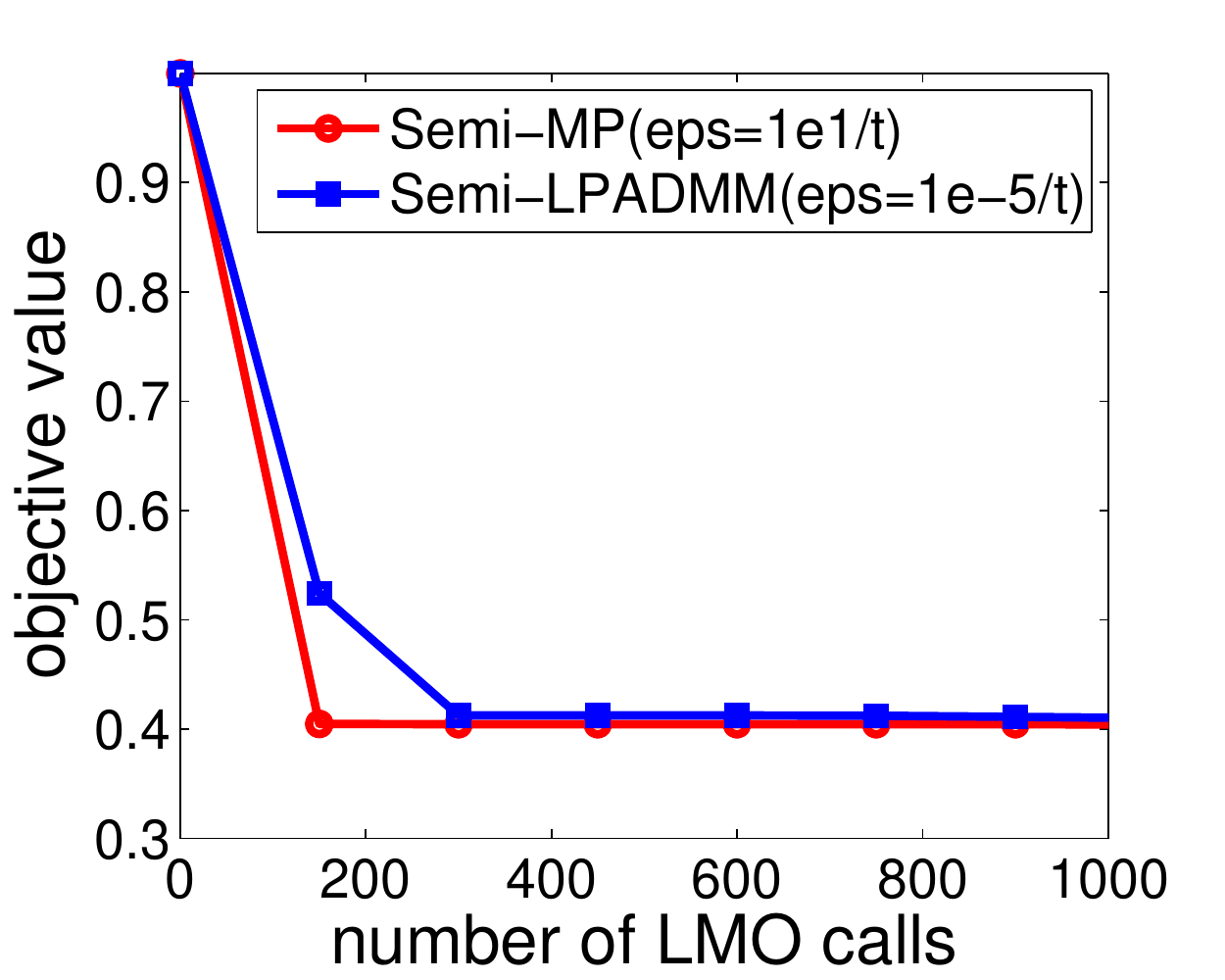}
  \end{minipage}
\caption{Robust collaborative filtering and link prediction: objective function vs elapsed time.\\  From left to right: (a) MovieLens 100K; (b) MovieLens 1M; (c) Wikivote(1024); (d) Wikivote(full)}
\label{fig:all-figs}
\end{figure}

\paragraph{Robust collaborative filtering}
We consider the collaborative filtering problem, with a nuclear-norm regularisation penalty and an $\ell_1$-loss function. 
We run the above three algorithms on the the small and medium MovieLens datasets. The small-size dataset consists of 943 users and 1682 movies with about 100K ratings, while the medium-size dataset consists of 3952 users and 6040 movies with about 1M ratings. We follow~\cite{pierucci2014smoothing} to set the regularisation parameters. 
 In Fig.~\ref{fig:all-figs}, we can see that Semi-MP clearly outperforms Smooth-CG, while it is competitive with Semi-SPG. 

\paragraph{Link prediction}
We consider now the link prediction problem, where the objective consists a hinge-loss for the empirical risk part and
multiple regularization penalties, namely the $\ell_1$-norm and the nuclear-norm. 
For this example, applying the Smooth-CG or Semi-SPG would require two smooth approximations, one for hinge loss term and one for $\ell_1$ norm term. Therefore, we consider another alternative approach, Semi-LPADMM, where we apply the linearized preconditioned ADMM algorithm~\cite{Lan2014} by solving proximal mapping through conditional gradient routines. Up to our knowledge, ADMM with early stopping is not fully theoretically analysed in literature. However, from an intuitive point of view, as long as the accumulated error is controlled sufficiently, such variant of ADMM should converge. 

We conduct experiments on a binary social graph data set called Wikivote, which consists of 7118 nodes and 103,747 edges. Since the computation cost of these two algorithms mainly come from the LMO calls, we present in below the performance in terms of number of LMO calls. For the first set of experiments, we select top 1024 highest degree users from Wikivote and run the two algorithms on this small dataset with different strategies for the inner LMO calls. 

In Fig.~\ref{fig:all-figs}, we observe that the Semi-MP is less sensitive to the inner accuracies of prox-mappings compared to the ADMM variant, which sometimes stops progressing if the prox-mapping of early iterations are not solved with sufficient accuracy.  
The results on the full dataset corroborate the fact that Semi-MP outperforms the semi-proximal variant of the ADMM algorithm.

\bibliographystyle{plain}

\begin{thebibliography}{10}

\bibitem{b15}
Francis Bach.
\newblock Duality between subgradient and conditional gradient methods.
\newblock {\em SIAM Journal on Optimization}, 2015.

\bibitem{Bach12}
Francis Bach, Rodolphe Jenatton, Julien Mairal, and Guillaume Obozinski.
\newblock Optimization with sparsity-inducing penalties.
\newblock {\em Found. Trends Mach. Learn.}, 4(1):1--106, 2012.

\bibitem{Bauschke:2011}
H.~H. Bauschke and P.~L. Combettes.
\newblock {\em Convex Analysis and Monotone Operator Theory in Hilbert Spaces}.
\newblock Springer, 2011.

\bibitem{bertsekas:2015}
D.~P. Bertsekas.
\newblock {\em Convex Optimization Algorithms}.
\newblock Athena Scientific, 2015.

\bibitem{byrd1995limited}
Richard~H Byrd, Peihuang Lu, Jorge Nocedal, and Ciyou Zhu.
\newblock A limited memory algorithm for bound constrained optimization.
\newblock {\em SIAM Journal on Scientific Computing}, 16(5):1190--1208, 1995.

\bibitem{Chen12}
Xi~Chen, Qihang Lin, Seyoung Kim, Jaime~G Carbonell, and Eric~P Xing.
\newblock Smoothing proximal gradient method for general structured sparse
  regression.
\newblock {\em The Annals of Applied Statistics}, 6(2):719--752, 2012.

\bibitem{cox2013dual}
Bruce Cox, Anatoli Juditsky, and Arkadi Nemirovski.
\newblock Dual subgradient algorithms for large-scale nonsmooth learning
  problems.
\newblock {\em Mathematical Programming}, pages 1--38, 2013.

\bibitem{dudik:harchaoui:malik2012}
M.~Dudik, Z.~Harchaoui, and J.~Malick.
\newblock Lifted coordinate descent for learning with trace-norm
  regularization.
\newblock {\em {Proceedings of the 15th International Conference on Artificial
  Intelligence and Statistics (AISTATS)}}, 2012.

\bibitem{garber2013linearly}
Dan Garber and Elad Hazan.
\newblock A linearly convergent conditional gradient algorithm with
  applications to online and stochastic optimization.
\newblock {\em arXiv preprint arXiv:1301.4666}, 2013.

\bibitem{harchaoui2013conditional}
Zaid Harchaoui, Anatoli Juditsky, and Arkadi Nemirovski.
\newblock Conditional gradient algorithms for norm-regularized smooth convex
  optimization.
\newblock {\em Mathematical Programming}, pages 1--38, 2013.

\bibitem{hazan2012projection}
E.~Hazan and S.~Kale.
\newblock Projection-free online learning.
\newblock In {\em ICML}, 2012.

\bibitem{CoMP13}
Niao He, Anatoli Juditsky, and Arkadi Nemirovski.
\newblock Mirror prox algorithm for multi-term composite minimization and
  semi-separable problems.
\newblock {\em arXiv preprint arXiv:1311.1098}, 2013.

\bibitem{J13}
Martin Jaggi.
\newblock Revisiting {F}rank-{W}olfe: Projection-free sparse convex
  optimization.
\newblock In {\em ICML}, pages 427--435, 2013.

\bibitem{juditsky13variational}
Anatoli Juditsky and Arkadi Nemirovski.
\newblock Solving variational inequalities with monotone operators on domains
  given by linear minimization oracles.
\newblock {\em arXiv preprint arXiv:1312.107}, 2013.

\bibitem{lan2013complexity}
Guanghui Lan.
\newblock The complexity of large-scale convex programming under a linear
  optimization oracle.
\newblock {\em arXiv}, 2013.

\bibitem{Lan14}
Guanghui Lan and Yi~Zhou.
\newblock Conditional gradient sliding for convex optimization.
\newblock {\em arXiv}, 2014.

\bibitem{mu2014scalable}
Cun Mu, Yuqian Zhang, John Wright, and Donald Goldfarb.
\newblock Scalable robust matrix recovery: Frank-wolfe meets proximal methods.
\newblock {\em arXiv preprint arXiv:1403.7588}, 2014.

\bibitem{Nem04}
Arkadi Nemirovski.
\newblock Prox-method with rate of convergence $o (1/t)$ for variational
  inequalities with lipschitz continuous monotone operators and smooth
  convex-concave saddle point problems.
\newblock {\em SIAM Journal on Optimization}, 15(1):229--251, 2004.

\bibitem{Nem10}
Arkadi Nemirovski, Shmuel Onn, and Uriel~G Rothblum.
\newblock Accuracy certificates for computational problems with convex
  structure.
\newblock {\em Mathematics of Operations Research}, 35(1):52--78, 2010.

\bibitem{Nesterov07a}
Y.~Nesterov.
\newblock Smoothing technique and its applications in semidefinite
  optimization.
\newblock {\em Math. Program.}, 110(2):245--259, 2007.

\bibitem{nesterov2005smooth}
Yu~Nesterov.
\newblock Smooth minimization of non-smooth functions.
\newblock {\em Mathematical programming}, 103(1):127--152, 2005.

\bibitem{Lan2014}
Yuyuan Ouyang, Yunmei Chen, Guanghui Lan, and Eduardo Pasiliao~Jr.
\newblock An accelerated linearized alternating direction method of
  multipliers, 2014.
\newblock \url{http://arxiv.org/abs/1401.6607}.

\bibitem{Parikh13}
Neal Parikh and Stephen Boyd.
\newblock Proximal algorithms.
\newblock {\em Foundations and Trends in Optimization}, pages 1--96, 2013.

\bibitem{pierucci2014smoothing}
Federico Pierucci, Zaid Harchaoui, and J{\'e}r{\^o}me Malick.
\newblock A smoothing approach for composite conditional gradient with
  nonsmooth loss.
\newblock In {\em Conf{\'e}rence d’Apprentissage Automatique--Actes
  CAP’14}, 2014.

\bibitem{inexact-prox-grad}
Mark Schmidt, Nicolas~L. Roux, and Francis~R. Bach.
\newblock Convergence rates of inexact proximal-gradient methods for convex
  optimization.
\newblock In {\em Adv. NIPS}. 2011.

\bibitem{ZYS12}
X.~Zhang, Y.~Yu, and D.~Schuurmans.
\newblock Accelerated training for matrix-norm regularization: A boosting
  approach.
\newblock In {\em NIPS}, 2012.

\end{thebibliography}

%
\newpage
\appendix
In this Appendix, we provide additional material on variational inequalities
and non-smooth optimisation algorithms, give the proofs on the main theorems, and provide additional information 
regarding the competing algorithms based on smoothing techniques and the implementation details for different models.

\section{Preliminaries: Variational Inequalities and Accuracy Certificates}\label{sect:preliminaries}
For the reader's convenience, we recall here the relationship between variational inequalities, 
accuracy certificates, and execution protocols, for non-smooth optimization algorithms. The exposition
below is directly taken from~\cite{CoMP13}, and recalled here for the reader's convenience. 

\paragraph{Execution protocols and accuracy certificates.}
Let $X$ be a nonempty closed convex set in a Euclidean space $E$ and $F(x):X\to E$ be a vector field.
\par
Suppose that we process $(X,F)$ by an algorithm which generates
a sequence of search points $x_t\in X$, $t=1,2,...$, and computes the vectors $F(x_t)$, so that after $t$ steps we have at our disposal {\sl $t$-step execution protocol}
$\cI_t=\{x_\tau,F(x_\tau)\}_{\tau=1}^t$. By definition, an {\sl accuracy certificate} for this protocol is simply a collection $\lambda^t=\{\lambda^t_\tau\}_{\tau=1}^t$ of nonnegative reals summing up to 1.
We associate with the protocol $\cI_t$ and accuracy certificate $\lambda^t$ two quantities as follows:
\begin{itemize}
\item {\sl Approximate solution} $x^t(\cI_t,\lambda^t):=\sum_{\tau=1}^t \lambda^t_\tau x_\tau$, which is a point of $X$;
\item {\sl Resolution $\Res(X'\big|\cI_t,\lambda^t)$ on a subset $X'\neq\emptyset$} of $X$ given by
\begin{equation}\label{resolution}
\Res(X'\big|\cI_t,\lambda^t) = \sup\limits_{x\in X'}\sum_{\tau=1}^t\lambda^t_\tau\langle F(x_\tau),x_\tau-x\rangle.
\end{equation}
\end{itemize}
The role of those notions for non-smooth optimization is explained below. 

\paragraph{Variational inequalities.} Assume that $F$ is {\em monotone}, i.e.,\text{VI}(X,F)
\begin{equation}\label{eq:near-monotone}
\langle F(x)-F(y),x-y\rangle\ge 0, \;\;\forall x,y\in X \; .
\end{equation}
Our goal is to approximate a weak solution
to the variational inequality (v.i.) $\text{VI}(X,F)$ associated with $(X,F)$. A weak solution is defined as a point $x_*\in X$ such that
\begin{equation}\label{vin}
\langle F(y),y-x_*\rangle \geq0\,\,\forall y\in X.
\end{equation}
A natural (in)accuracy measure  of a candidate weak solution $x\in X$ to $\text{VI}(X,F)$ is the {\sl dual gap function}
\begin{equation}\label{dualgap}
\epsilonvi(x\big|X,F) = \sup_{y\in X} \langle F(y),x-y\rangle
\end{equation}
This inaccuracy is a convex nonnegative function which vanishes exactly at the set of weak solutions to the $\text{VI}(X,F)$.
\begin{proposition}\label{prop1} For every $t$, every execution protocol $\cI_t=\{x_\tau\in X,F(x_\tau)\}_{\tau=1}^t$ and every
accuracy certificate $\lambda^t$ one has $x^t:=x^t(\cI_t,\lambda^t)\in X$. Besides this, {assuming $F$ monotone,} for  every closed convex set $X'\subset X$ such that $x^t\in X'$ one has
\begin{equation}\label{epsilonvi}
\epsilonvi(x^t\big|X',F)\leq \Res(X'\big|\cI_t,\lambda^t).
\end{equation}
\end{proposition}
{\bf Proof.} Indeed, $x^t$ is a convex combination of the points $x_\tau\in X$ with coefficients $\lambda^t_\tau$, whence $x^t\in X$. With $X'$ as in the premise of Proposition, we have
$$
\forall y\in X': \langle F(y),x^t-y\rangle =\sum_{\tau=1}^t\lambda^t_\tau\langle F(y),x_\tau-y\rangle \leq \sum_{\tau=1}^t \lambda^t_\tau \langle F(x_\tau),x_\tau-y\rangle
\leq\Res(X'\big|\cI_t,\lambda^t),
$$
where the {first $\leq$ is due to} monotonicity of $F$. \qed
\paragraph{Convex-concave saddle point problems.} Now let $X=X_1\times X_2$, where $X_i$ is a closed convex subset in Euclidean space $E_i$, $i=1,2$, and $E=E_1\times E_2$, and let
$\Phi(x^1,x^2):X_1\times X_2\to\bR$
be a locally Lipschitz continuous function which is convex in $x^1\in X_1$ and concave in $x^2\in X_2$. $X_1,X_2,\Phi$ give rise to  the saddle point problem
\begin{equation}\label{saddlepoint}
\SadVal=\min_{x^1\in X_1}\max_{x^2\in X_2}\Phi(x^1,x^2),
\end{equation} two induced convex optimization problems
\begin{equation}\label{induced}
\begin{array}{rclr}
\Opt(P)&=&\displaystyle\min_{x^1\in X_1}\bigg[\overline{\Phi}(x^1)=\sup_{x^2\in X_2}\Phi(x^1,x^2)\bigg]&(P)\\
\Opt(D)&=&\displaystyle\max_{x^2\in X_2}\bigg[\underline{\Phi}(x^2)=\inf_{x^1\in X_1}\Phi(x^1,x^2)\bigg]&(D)\\
\end{array}
\end{equation}
and a vector field $F(x^1,x^2)=[F_1(x^1,x^2);F_2(x^1,x^2)]$ specified (in general, non-uniquely) by the relations
$$
\forall (x^1,x^2)\in X_1\times X_2: F_1(x^1,x^2)\in\partial_{x^1}\Phi(x^1,x^2),\,F_2(x^1,x^2)\in\partial_{x^2}[-\Phi(x^1,x^2)].
$$
It is well known that $F$ is monotone on $X$, and that weak solutions to the $\text{VI}(X,F)$ are exactly the saddle points of $\Phi$ on $X_1\times X_2$. These saddle points exist if and only if
$(P)$ and $(D)$ are solvable with equal optimal values, in which case the saddle points are exactly the pairs $(x^1_*,x^2_*)$ comprised by optimal solutions to $(P)$ and $(D)$.
In general, $\Opt(P)\geq\Opt(D)$, with equality definitely taking place when at least one of the sets $X_1,X_2$ is bounded; if both are bounded, saddle points do exist. {To avoid unnecessary complications, from now on, when speaking about a convex-concave saddle point problem, we assume that the problem is {\sl proper}, meaning that $\Opt(P)$ and $\Opt(D)$ are reals; this definitely is the case when $X$ is bounded.
\par
A natural (in)accuracy measure for a candidate $x=[x^1;x^2]\in X_1\times X_2$ to the role of a saddle point of $\Phi$ is the quantity
\begin{equation}\label{saddlepointinacc}
\begin{array}{rcl}
\epsilonsad(x\big|X_1,X_2,\Phi)&=&\overline{\Phi}(x^1) -\underline{\Phi}(x^2)\\
&=&[\overline{\Phi}(x^1)-\Opt(P)]+[\Opt(D)-\underline{\Phi}(x^2)] +\underbrace{[\Opt(P)-\Opt(D)]}_{\geq0}\\
\end{array}
\end{equation}
This inaccuracy is nonnegative and is the sum of the duality gap $\Opt(P)-\Opt(D)$ (always nonnegative and vanishing when one of the sets $X_1,X_2$ is bounded)  and the inaccuracies,
in terms of respective objectives, of $x^1$ as a candidate solution to $(P)$ and $x^2$ as a candidate solution to $(D)$.
\par
The role of accuracy certificates in convex-concave saddle point problems stems from the following observation:
\begin{proposition} \label{prop2} Let $X_1,X_2$ be nonempty closed convex sets, $\Phi:X:=X_1\times X_2\to\bR$ be a locally Lipschitz continuous convex-concave function, and
$F$ be the associated monotone vector field on $X$.
\par Let $\cI_t=\{x_\tau=[x^1_\tau;x^2_\tau]\in X,{F}(x_\tau)\}_{\tau=1}^t$ be a $t$-step execution protocol associated with $(X,{F})$ and
$\lambda^t=\{\lambda^t_\tau\}_{\tau=1}^t$ be an associated
accuracy certificate. Then $x^t:=x^t(\cI_t,\lambda^t)=[x^{1,t};x^{2,t}] \in X$.
\par
Assume, further, that $X^\prime_1\subset X_1$ and $X^\prime_2\subset X_2$ are closed convex sets such that
\begin{equation}\label{condition}
x^t\in X^\prime:=X^\prime_1\times X^\prime_2.
\end{equation}
Then
\begin{equation}\label{saddlepointres}
\epsilonsad(x^t\big|X_1^\prime,X_2^\prime,\Phi)=\sup_{x^2\in X_2'}\Phi(x^{1,t},x^2)-\inf_{x^1\in X_1'}\Phi(x^1,x^{2,t})\leq\Res(X^\prime\big|\cI_t,\lambda^t).
\end{equation}
In addition, setting $\widetilde{\Phi}(x^1)=\sup_{x^2\in X_2^\prime}\Phi(x^1,x^2)$, for every $\bar{x}^1\in X^\prime_1$ we have
\begin{equation}\label{forfeasible}
\widetilde{\Phi}(x^{1,t})-\widetilde{\Phi}(\bar{x}^1)\leq\widetilde{\Phi}(x^{1,t})-\Phi(\bar{x}^1,x^{2,t})\leq\Res(\{\bar{x}^1\}\times X^\prime_2\big|\cI_t,\lambda^t).
\end{equation}
In particular, when the problem $\Opt=\min_{x^1\in X^\prime_1}\widetilde{\Phi}(x^1)$ is solvable with an optimal solution $x^1_*$, we have
\begin{equation}\label{ifsolvable}
\widetilde{\Phi}(x^{1,t})-\Opt\leq\Res(\{x^1_*\}\times X^\prime_2\big|\cI_t,\lambda^t).
\end{equation}
\end{proposition}
{\bf Proof.} The inclusion $x^t\in X$ is clear. For every set $Y\subset X$ we have
\begin{align*}
\begin{array}{l}
\forall [p;q]\in Y:\\
\Res(Y\big|\cI_t,\lambda^t)\geq \sum_{\tau=1}^t\lambda^t_\tau\left[\langle F_1(x^1_\tau),x^1_\tau-p\rangle + \langle F_2(x^2_\tau),x^2_\tau-q\rangle \right]\\
\geq \sum_{\tau=1}^t\lambda^t_\tau\left[[\Phi(x^1_\tau,x^2_\tau)-\Phi(p,x^2_\tau)]+[\Phi(x^1_\tau,q)-\Phi(x^1_\tau,x^2_\tau)]\right]\\
\hbox{[by the origin of $F$ and since $\Phi$ is convex-concave]}\\
=\sum_{\tau=1}^t\lambda^t_\tau\left[\Phi(x^1_\tau,q)-\Phi(p,x^2_\tau)\right]
\geq \Phi(x^{1,t},q)-\Phi(p,x^{2,t})\\
\;~\;~\;~\;~\;~\;~\;~\;~\;~\hbox{[by origin of $x^t$ and since $\Phi$ is convex-concave]}\\
\end{array}
\end{align*}
Thus, for every $Y\subset X$ we have
\begin{equation}\label{conclude}
\sup_{[p;q]\in Y} \left[\Phi(x^{1,t},q)-\Phi(p,x^{2,t})\right]\leq \Res(Y\big|\cI_t,\lambda^t).
\end{equation}
Now assume that Condition (\ref{condition}) is satisfied. Setting $Y=X':=X^\prime_1\times X^\prime_2$, and recalling what $\epsilonsad$ is, (\ref{conclude}) yields (\ref{saddlepointres}). With $Y=\{\bar{x}^1\}\times X^\prime_2$ (\ref{conclude}) yields the second inequality in
(\ref{forfeasible}); the first inequality in (\ref{forfeasible}) is clear since $x^{2,t}\in X^\prime_2$. \qed

\section{Theoretical analysis of composite Mirror Prox with inexact proximal mappings}
We restate the Theorem~\ref{thm:theMP} below and the proof below. The theoretical convergence rate 
established in Theorem~\ref{thm:theMP} and Corollary~\ref{cor:theMP} extends the previous result established in Corollary 3.1 in~\cite{CoMP13} for CMP with exact prox-mappings. Indeed, when exact prox-mappings are used, we recover the result of~\cite{CoMP13}. When inexact prox-mappings are used, 
the errors due to the inexactness of the prox-mappings accumulates and is reflected in the bound (\ref{epsilonvismall}) and (\ref{optimalitygap}). 

\begingroup
\def\thetheorem{\ref{thm:theMP}}
\begin{theorem}
Assume that the sequence of step-sizes $(\gamma_t)$ in the CMP algorithm satisfy
\begin{equation}\label{gammaupperbound}
\sigma_t:=\gamma_t\langle F_u(\hat{u}^t)-F_u(u^t),\hat{u}^t-u^{t+1}\rangle-V_{\hat{u}^t}(u^{t+1})-V_{u^t}(\hat{u}^{t})\le \gamma_t^2M^2 \, , \quad t=1,2,\ldots, T \; .
\end{equation}
Then, denoting $\Theta[X]=\sup_{[u;v]\in X}V_{u^1}(u)$, for a sequence of inexact prox-mappings with inexactness $\epsilon_t\geq 0$, we have
\begin{equation}\label{epsilonvismall}
\epsilonvi(\bar{x}_T\big|X,F):=\sup_{x\in X} \; \left\langle F(x), \bar{x}_T-x\right\rangle \leq \frac{\Theta[X]+M^2{\sum}_{t=1}^T\gamma_t^2+2{\sum}_{t=1}^T\epsilon_t}{\sum_{t=1}^T\gamma_t} \; .
\end{equation}
\end{theorem}
\addtocounter{theorem}{-1}
\endgroup

\textbf{Remarks }
Note that the assumption on the sequence of step-sizes $(\gamma_t)$ is clearly satisfied when $\gamma_t\leq ({\sqrt{2}L})^{-1}$. When $M=0$, it is satisfied as long as $\gamma_t\leq L^{-1}$. 

\begin{proof}
The proofs builds upon and extends the proof in~\cite{CoMP13}. For all $u,u',w\in U$, we have the well-known identity
\begin{equation}\label{threetermid}
\langle V'_{u}(u'),w-u'\rangle =V_{u}(w)-V_{u'}(w)-V_{u}(u').
\end{equation}

Indeed, the right hand side writes as
\begin{align*}
\lefteqn{[\omega(w)-\omega(u)-\langle\omega'(u),w-u\rangle]-[\omega(w)-\omega(u')-\langle\omega'(u'),w-u'\rangle]-[\omega(u')-\omega(u)-\langle\omega'(u),u'-u\rangle]}\\
&=&\langle \omega'(u),u-w\rangle +\langle \omega'(u),u'-u\rangle +\langle \omega'(u'),w-u'\rangle =\langle\omega'(u')- \omega'(u),w-u'\rangle=\langle V'_{u}(u'),w-u'\rangle.
\end{align*}
For $x=[u;v]\in X,\;\xi=[\eta;\zeta]$, $\epsilon\geq0$,  let $[u';v']\in P_x^\epsilon(\xi)$. By definition, for all $[s;w]\in X$, the inequality holds
\[
\langle \eta+V'_u(u'),u'-s\rangle+\langle\zeta,v'-w\rangle \le \epsilon,
\]
which by \rf{threetermid} implies that
\begin{equation}\label{prox_lemma}
\langle \eta,u'-s\rangle+\langle\zeta,v'-w\rangle \le \langle V'_u(u'),s-u'\rangle+\epsilon=
V_{u}(s)-V_{u'}(s)-V_{u}(u')+\epsilon.
\end{equation}
 When applying \rf{prox_lemma} with $\epsilon=\epsilon_t$, $[u;v]=[u^{t};v^{t}]=x^{t}$, $\xi=\gamma_{t} F(x^{t})=[\gamma_{t} F_u(u^{t});\gamma_{t} F_v]$, $[{u}';{v}']=[\hat{u}^{t};\hat v^{t}]=y^{t}$, and $[s;w]=[u^{{t}+1};v^{{t}+1}]=x^{{t}+1}$ we obtain
\begin{equation}
\label{prox100}
\gamma_{t} [\langle F_u(u^{t}),\hat{u}^{t}-u^{{t}+1}\rangle+\langle F_v,\hat{v}^{t}-v^{{t}+1}\rangle]\le V_{u^{t}}(u^{{t}+1})-V_{\hat{u}^{t}}(u^{{t}+1})-V_{u^{t}}(\hat{u}^{{t}})+\epsilon_{t} \; ; 
\end{equation}
and applying \rf{prox_lemma} with $\epsilon=\epsilon_{t}$, $[u;v]=x^{t}$, $\xi=\gamma_{t} F(y^{t})$, $[{u}';{v}']=x^{{t}+1}$, and $[s;w]=z\in X$ we get
\begin{equation}
\label{prox101}
\gamma_{t} [\langle F_u(\hat{u}^{t}),u^{{t}+1}-s\rangle+\langle F_v,v^{{t}+1}-w\rangle]\le V_{u^{t}}(s)-V_{u^{{t}+1}}(s)-V_{u^{t}}(u^{{t}+1})+\epsilon_{t} \; .
\end{equation}
Adding \rf{prox101} to \rf{prox100}, we obtain for every $z=[s;w]\in X$
\begin{align}
\label{prox102}
\gamma_{t} \langle F(y^{t}),y^{t}-z\rangle &=\gamma_{t} [\langle F_u(\hat{u}^{t}),\hat{u}^{t}-s\rangle+\langle F_v,\hat{v}^{t}-w\rangle]  \nonumber \\
&\le V_{u^{t}}(s)-V_{u^{{t}+1}}(s)+\sigma_{t}+2\epsilon_{t} \; ,
\end{align}
with 
\begin{equation*}
\sigma_{t} := \gamma_{t}\langle F_u(\hat{u}^{t})-F_u(u^{t}),\hat{u}^{t}-u^{{t}+1}\rangle-V_{\hat{u}^{t}}(u^{{t}+1})-V_{u^{t}}(\hat{u}_{{t}}) \; .
\end{equation*}
Due to the strong convexity, with modulus 1, of $V_u(\cdot)$ w.r.t. $\|\cdot\|$, we have for all $u,\hat{u}$
\begin{equation*}
V_u(\hat{u})\geq {1\over 2}\|u-\hat{u}\|^2 \; .
\end{equation*} 
Therefore,
\bse
\sigma_{t}&\leq& \gamma_{t}\|F_u(\hat{u}^{t})-F_u(u^{t})\|_*\|\hat{u}^{t}-u^{{t}+1}\|-\half\|\hat{u}^{t}-u^{{t}+1}\|^2-\half\|u^{t}-\hat{u}^{t}\|^2\\
&\leq&
\half\left[\gamma_{t}^2\|F_u(\hat{u}^{t})-F_u(u^{t})\|_*^2-\|u^{t}-\hat{u}^{t}\|^2\right]\\
&\leq&
\half\left[\gamma_{t}^2[M+L\|\hat{u}^{t}-u^{t}\|]^2-\|u^{t}-\hat{u}^{t}\|^2\right],
\ese
where the last inequality follows from Assumption \textbf{A}.3. Note that $\gamma_{t} L<1$ implies that
\begin{equation*}
\gamma_{t}^2[M+L\|\hat{u}^{t}-u^{t}\|]^2-\|\hat{u}^{t}-u^{t}\|^2\leq \max_r
\left[\gamma_{t}^2[M+Lr]^2-r^2\right]={\gamma_{t}^2M^2\over1-\gamma_{t}^2L^2}.
\end{equation*} 
Let us assume that the step-sizes $\gamma_{t}>0$ are chosen so that \rf{gammaupperbound} holds, that is 
$\sigma_{t} \leq \gamma_{t}^2M^2$. 
It is indeed the case when $0<\gamma_{t}\leq {1\over\sqrt{2}L}$; when $M=0$, we can take also $\gamma_{t}\leq {1\over L}$.
Summing up inequalities~(\ref{prox102}) over ${t}=1,2,...,t$, and taking into account that $V_{u^{t+1}}(s)\geq0$, we finally conclude that for all $z=[s;w]\in X$,
\begin{align*}
\sum_{{t}=1}^T\lambda_T^{t}\langle F(y^{t}),y^{t}-z\rangle 
&\leq{V_{u^1}(s) +M^2\sum_{{t}=1}^T\gamma_{t}^2+2\sum_{t=1}^T\epsilon_t\over
\sum_{{t}=1}^T\gamma_{t}}, \text{ where }\lambda_T^{t}=(\sum_{i=1}^T\gamma_i )^{-1}\gamma_{t} \; .
\end{align*}
\end{proof}

\section{Theoretical analysis of composite conditional gradient}\label{sect:prooftheCCG}

\subsection{Convergence rate}
The CCG algorithm enjoys a convergence rate in $O(t^{-(\kappa-1)})$ in the evaluations of the function $\phi^+$, 
and the accuracy certificates $(\delta_t)$ enjoy the same rate $O(t^{-(\kappa-1)})$ as well, for solving problems of type~(\ref{CCGProblem}).
\begingroup
\def\theproposition{\ref{prop:theCCG}}
\begin{proposition}
Denote $D$ the $\|\cdot\|$-diameter of $U$. When solving problems of type~(\ref{CCGProblem}), the sequence of iterates $(x^t)$ of CCG satisfies
\begin{equation}\label{suboptimality}
\epsilon_t:=\phi^+(x^t)-\displaystyle\min_{x\in X}\phi^+(x)\leq {2L_0D^\kappa\over\kappa(3-\kappa)}\left(\frac{2}{t+1}\right)^{\kappa-1},\,t\geq2
\end{equation}
In addition, the accuracy certificates $(\delta_t)$ satisfy
\begin{equation}\label{certificate}
\min_{1\leq s\leq t} \; \delta_s \leq O(1)L_0D^\kappa \left(\frac{2}{t+1}\right)^{\kappa-1},\,t\geq2
\end{equation}
\end{proposition}
\endgroup

\subsection{Proof of Proposition \ref{prop:theCCG}}
\paragraph{1$^0$.} The projection of $X_2$ onto $E_{u_2}$ is contained in $U_2$, whence
$$
\|u_2[\nabla\phi(u_2^s)]-u_2^s\|\leq D.
$$
This observation, due to the structure of $\phi^+$, implies that whenever $x,x'\in X$ and $\gamma\in[0,1]$, we have
\begin{equation}\label{CCGWeHave}
\phi^+(x+\gamma(x^+-x)) \leq \phi^+(x)+\gamma\langle\nabla \phi^+(x),x'-x\rangle+{L_0D^\kappa\over\kappa}\gamma^\kappa.
\end{equation}
Setting $x^s_+=x_2^s+\gamma_s(x_2[\nabla\phi(u^s)]-x_2^s)$ and $\gamma_s 2/(s+1)$, we have
\begin{align}
\delta_{t+1}&\leq\phi^+(x^s_+)-\min_{x_2\in X_2}\phi^+(x_2) \\
&\leq\delta_s+\gamma_s\langle \nabla\phi(x_2^s),x[\nabla\phi^+(x_2^s)]-x_2\rangle +\frac{L_0D^\kappa}{\kappa}\gamma_s^\kappa
\\
&=\delta_s-\gamma_s\Delta^s+\frac{L_0D^\kappa}{\kappa}\gamma_s^\kappa,
\end{align}
whence, due to $\Delta_s\geq\delta_s\geq0$, 
\newcommand\numberthis{\addtocounter{equation}{1}\tag{\theequation}}
\begin{align*}
(i)\quad \delta_{t+1}&\leq(1-\gamma_s)\delta_s+\frac{L_0D^\kappa}{\kappa}\gamma_s^\kappa,\,s=1,2,...,\\
(ii)\quad \gamma_\tau\Delta_\tau&\leq\delta_\tau-\delta_{\tau+1}+\frac{L_0D^\kappa}{\kappa}\gamma_\tau^\kappa,\,\tau=1,2,...\numberthis\label{CCGWhence}
\end{align*}
\paragraph{2$^0$.} Let us prove (\ref{suboptimality}) by induction on $s\geq2$. 
By (\ref{CCGWhence}.i) and due to $\gamma_1=1$ we have $\delta_2\leq {L_0D^\kappa\over\kappa}$, whence $\delta_2\leq {2L_0D^\kappa\over\kappa(3-\kappa)}\gamma_2^{\kappa-1}$ due to $\gamma_2=2/3$ and $1<\kappa\leq 2$. Now assume that $\delta_s\leq {2L_0D^\kappa\over\kappa(3-\kappa)}\gamma_s^{\kappa-1}$ for some $t\geq2$. Then, invoking (\ref{CCGWhence}.i),
\begin{align*}
\delta_{s+1} &\leq {2L_0D^\kappa\over\kappa(3-\kappa)}\gamma_s^{\kappa-1}(1-\gamma_s)
+{L_0D^\kappa\over\kappa}\gamma_s^\kappa \\
&\leq {2L_0D^\kappa\over\kappa(3-\kappa)}\left[\gamma_s^{\kappa-1}-{\kappa-1\over 2}\gamma_s^\kappa\right] \\
&\leq{2L_0D^\kappa\over\kappa(3-\kappa)}2^{\kappa-1}{\left[(t+1)^{1-\kappa}+(1-\kappa)(t+1)^{-\kappa}
\right]}
\end{align*}
Therefore, by convexity of $(t+1)^{1-\kappa}$ in $t$
\begin{align*}
\delta_{s+1} &\leq {2L_0D^\kappa\over\kappa(3-\kappa)}2^{\kappa-1}(t+2)^{1-\kappa}={2L_0D^\kappa\over\kappa(3-\kappa)}\gamma_{t+1}^{\kappa-1}
\end{align*}
The induction is completed. 
\paragraph{3$^0$.} To prove (\ref{certificate}), given $s\geq2$, let $s_-=\hbox{Ceil}(\max[2,s/2])$. Summing up inequalities (\ref{CCGWhence}.ii) over $s_-\leq \tau\leq s$, we get
\begin{align*}
\left(\min_{\tau\leq s} \; \Delta_\tau\right) \;
{\sum}_{\tau=s_-}^s \gamma_\tau &\leq\sum_{\tau=s_-}^s\gamma_\tau \Delta_\tau \: \leq \delta_{s_-} -\delta_{s+1}+{L_0D^\kappa\over 2}{\sum}_{\tau=s_-}^s\gamma_\tau^\kappa \: \leq O(1)L_0D^\kappa\gamma_s^{\kappa-1}
\end{align*}
and $\sum_{\tau=s_-}^s \gamma_\tau\geq O(1)$, and (\ref{certificate}) follows. \qed

\section{\spmp}
\subsection{Theoretical analysis for \spmp}\label{sect:proofthecomplexity}
We first restate Proposition \ref{prop:thecomplexity} and provide the proof below. 
\begingroup
\def\theproposition{\ref{prop:thecomplexity}}
\begin{proposition}
Under the assumption $(\textbf{A}.1)-(\textbf{A}.4)$ and $(\textbf{S}.1)-(\textbf{S}.3)$ with $M=0$,  for the  outlined algorithm to return an $\epsilon$-solution to the variational inequality $VI(X,F)$, the total number of Mirror Prox steps required does not exceed $O\left(\frac{L\Theta[X]}{\epsilon}\right)$, 
and the total number of calls to the Linear Minimization Oracle does not exceed 
\begin{equation*}
\cN =O(1)\left(\frac{L_0L^\kappa D^\kappa}{\epsilon^\kappa}\right)^{\frac{1}{\kappa-1}}\Theta[X].
\end{equation*}
In particular, if we use Euclidean proximal setup on $U_2$ with $\omega_2(\cdot)=\frac{1}{2}\|x_2\|^2$, which leads to $\kappa =2$ and $L_0=1$, then the number of LMO calls does not exceed $\cN=O(1)\left(L^2D^2(\Theta[X_1]+D^2\right)/\epsilon^2$.
\end{proposition}
\endgroup
\begin{proof}
Let us fix $N$ as the number of Mirror prox steps, and since $M=0$, from Theorem \ref{thm:theMP}, the efficiency estimate of the variational inequality implies that 
$$\epsilonvi(\bar x^N|X,F)\leq \frac{L(\Theta[X]+2\sum_{t=1}^N\epsilon_t)}{N}.$$
Let us fix $\epsilon_t=\frac{2\Theta[X]}{N}$ for each $t=1,\ldots, N$, then from Proposition \ref{prop:theCCG}, it takes at most $s=O(1)(\frac{L_0D^\kappa N}{\Theta[X]})^{1/(\kappa-1)}$ LMO oracles to generate a point such that $\Delta_s\leq \epsilon_t$. Moreover, we have
$$\epsilonvi(\bar x^N|X,F)\leq 2\frac{L\Theta[X]}{N}.$$ Therefore, to ensure $\epsilonvi(\bar x^N|X,F)\leq \epsilon$ for a given accuracy $\epsilon>0$, the number of Mirror Prox steps $N$ is at most $O(\frac{L\Theta[X]}{\epsilon})$ and the number of LMO calls on $X_2$ needed is at most 
$$\cN=O(1)\Big(\frac{L_0D^\kappa N}{\Theta[X]}\Big)^{1/(\kappa-1)}\cdot N=O(1)\Big(\frac{L_0L^\kappa D^\kappa}{\epsilon^\kappa}\Big)^{1/(\kappa-1)}\Theta[X].$$
In particular, if $\kappa=2$ and $L_0=1$, this quantity can be reduced to 
$$\cN=O(1)\frac{L^2D^2\Theta[X]}{\epsilon^2}.$$
\end{proof}

\subsection{Discussion of \spmp}
\label{sec:discut}
The proposed \spmp algorithm enjoys the \emph{optimal complexity bounds}, 
i.e. $O(1/\epsilon^2)$, in the number of calls to linear minimization oracle. Furthermore, \spmp generalizes
previously proposed approaches and improves upon them in special cases of problem (\ref{nonsmoothproblem}).

When there is no regularisation penalty, \spmp is more general than previous algorithms for solving the corresponding constrained non-smooth optimisation problem.~\spmp does not require assumptions on “favorable geometry” of dual domains $Z$ or simplicity of $\psi(\cdot)$ in (\ref{sadpoint}). 
When the regularisation is simply a norm (with no operator in front of the argument),~\spmp is competitive with previously proposed approaches~\cite{Lan14, pierucci2014smoothing} based on smoothing techniques. 

When the regularisation penalty is non-trivial, \spmp is the first proximal-free or conditional-gradient-type optimization algorithm, up to our knowledge.

\section{Numerical experiments and implementation details} 
\label{sec:compet-details}

\subsection{Matrix completion: $\ell_2$-fit +nuclear norm} 
We first consider the the following type of matrix completion problem, 
\begin{equation}\label{matrixcompletionL2}
\min _{x\in\bR^{m\times n}} \|{P_\Omega}x-b\|_2 +\lambda\|x\|_\nuc
\end{equation}
where $\|\cdot\|_\nuc$ stands for the nuclear norm and $P_\Omega x$ is the restriction of $x$ onto the cells $\Omega$.

\paragraph{Competing algorithms.} We compare the following three candidate algorithms, i) \spmp (\textbf{Semi-MP}) ; ii) conditional gradient after smoothing (\textbf{Smooth-CG}); iii) inexact accelerate proximal gradient after smoothing (\textbf{Semi-SPG}). We provide below the key steps of each algorithms. 
\begin{enumerate}
\item \textbf{Semi-MP}: this is shorted for our \spmp algorithm, we solve the saddle point reformulation given by
\begin{equation}\label{saddlepointL2}
\min _{x,v: \|x\|_\nuc\leq v}\max_{\|y\|_2\leq 1} \langle P_\Omega x-b, y\rangle +\lambda v
\end{equation}
which is equivalent as to the semi-structured variational inequality \semivi$(X,F)$ with 
$X=\{[u=(x;y);v]: \|x\|_\nuc\leq v,\|y\|_2\leq 1\}$ and $F=[F_u(u);F_v]=[P_\Omega^Ty; b-P_\Omega x;\lambda]$. The subdomain $X_1=\{y:\|y\|_2\leq 1\}$ is given by full-prox setup and  the subdomain $X_2=\{(x;v):\|x\|_\nuc\leq v\}$ is given by LMO. By setting both the distance generating functions $\omega_x(x)$ and $\omega_y(y)$ as the Euclidean distance, the update of $y$ reduces to a gradient step, and the update of $x$ follows the composite conditional gradient routine over a simple quadratic problem.

\item \textbf{Smooth-CG}: The algorithm (\cite{pierucci2014smoothing}) directly 
applies the generalized composite conditional gradient on the following smoothed problem using the Nesterov smoothing technique,
\begin{equation}\label{smoothL2}
\min _{x,v: \|x\|_\nuc\leq v}f^\gamma(x)+\lambda v, \text{ where } f^{\gamma}(x)=\max_{\|y\|_2\leq 1}\{\langle P_\Omega x-b,y\rangle-\frac{\gamma}{2}\|y\|_2^2\} .
\end{equation}

Under the full memory version, the update of $x$ at step $t$ requires computing reoptimization problem
\begin{equation}\label{reoptimize}
\min_{\theta_1,\ldots,\theta_t} f^{\gamma}(\sum_{i=1}^t\theta_i u_iv_i^T)+\lambda\sum_{i=1}^t\theta_i
\end{equation}
where $\{u_i,v_i\}_{i=1}^t$ are the singular vectors collected from the linear minimization oracles. Same as suggested in \cite{pierucci2014smoothing}, we use the quasi-Newton solver L-BFGS-B \cite{byrd1995limited} to solve the above re-optimization subproblem.  Notice that in this situation, solving (\ref{reoptimize}) can be relatively efficient even for large $t$ since computing the gradient of the objective in (\ref{reoptimize})   does not necessarily need to compute out the full matrix representation of $x=\sum_{i=1}^t\theta_iu_iv_i^T$.

\item \textbf{Semi-SPG}: The approach is to apply the accelerated proximal gradient to the smoothed composite model as in (\ref{smoothL2}) and approximately solve the proximal mappings via conditional gradient routines. In fact, Semi-SPG can be considered as a direct extension of the conditional gradient sliding
to the composite setting. Same as Semi-MP, the update of $x$ is given by the composite conditional gradient routine over a simple quadratic problem and additional  interpolation step. Since the Lipschitz constant is not known,  the learning rate is selected through backtracking. 
\end{enumerate}

For Semi-MP and Semi-SPG, we test two different strategies for the inexact prox-mappings, a)fixed inner CG steps and b)decaying $\epsilon_t=c/t$ as the theory suggested. For the sake of simplicity, we generate the synthetic data such that the magnitudes of the constant factors (i.e. Frobenius norm and nuclear norm of optimal solution) are approximately of order 1, which means the convergence rate is dominated mainly by the number of LMO calls.  In Fig.~\ref{fig:synthetic}, we evaluate the optimality gap of these algorithms with different parameters (e.g. number of inner steps, scaling factor c, smoothness parameter $\gamma$) and compare their performance given the best-tuned parameter. As the plot shows, the Semi-MP algorithm generates a solution with $\epsilon=10^{-3}$ accuracy within about 3000 LMO calls, which is not bad at all given the fact that the worst complexity is $O(1/\epsilon^2)$. Also, the plots indicate that using the second strategy with $O(1/t)$ decaying inexactness provides better and more reliable performance than using fixed number of inner steps. Similar trends are observed for the Semi-SPG. One can see that these two algorithms based on inexact proximal mappings are notably faster than applying conditional gradient on the smoothed problem. Moreover, since the Smooth-CG requires additional computation and memory cost for the re-optimization procedure, the actual difference in terms of CPU time could be more significant. 

\begin{figure} [!ht]
  \begin{minipage}[t]{.25\textwidth}
    \includegraphics[scale=0.3]{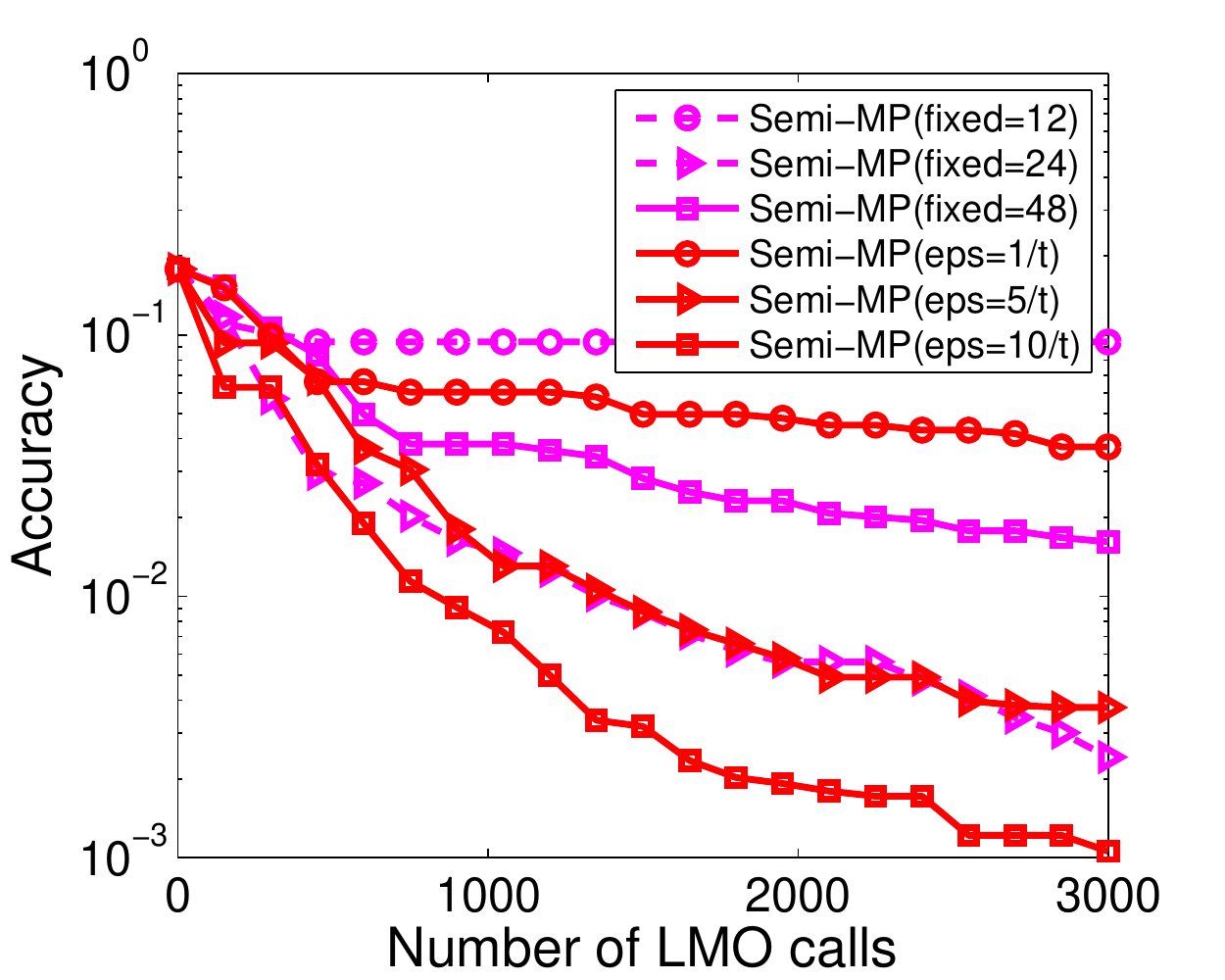}
  \end{minipage}%
  \begin{minipage}[t]{.25\textwidth}
    \includegraphics[scale=0.3]{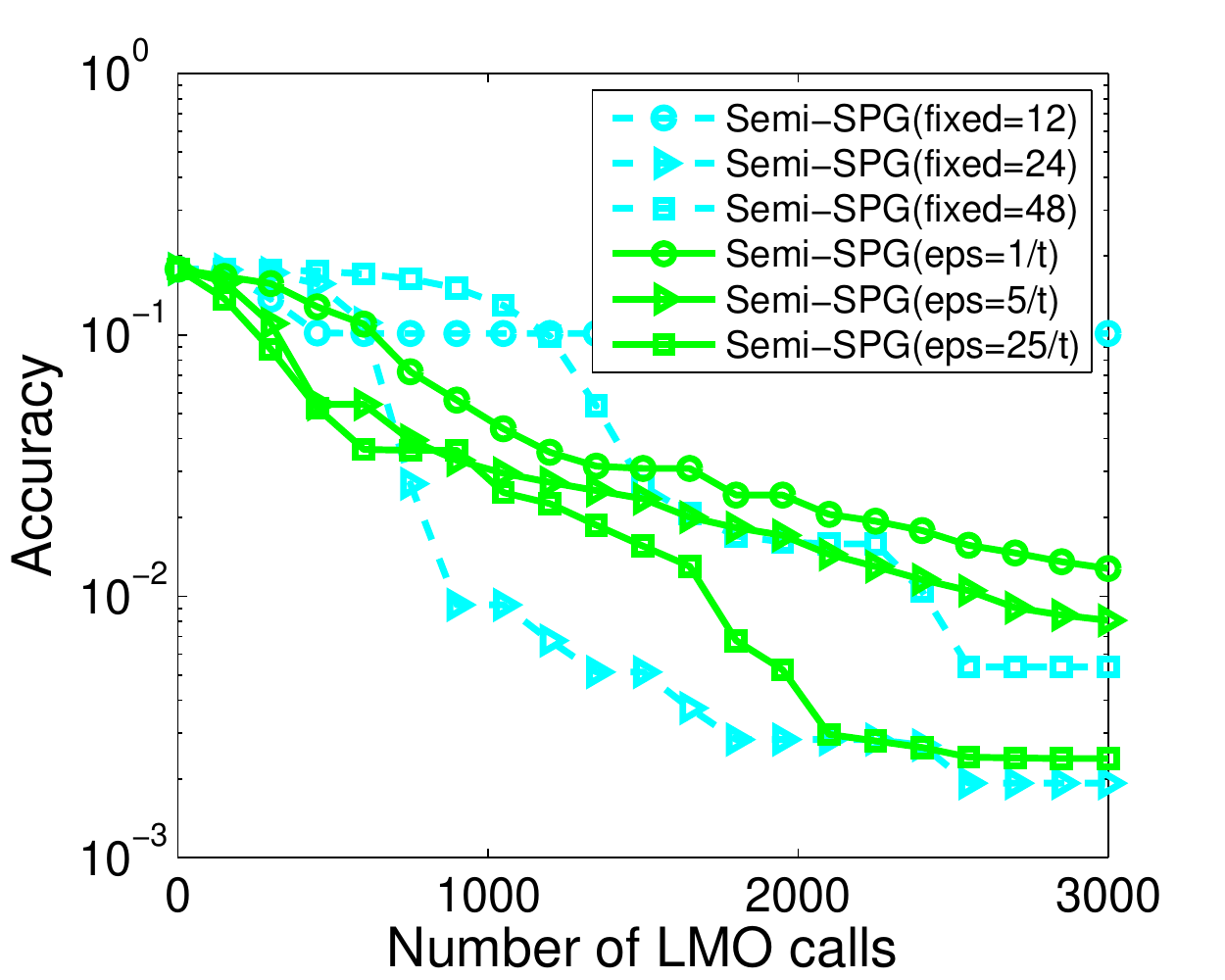}
  \end{minipage}%
  \begin{minipage}[t]{.25\textwidth}
    \includegraphics[scale=0.3]{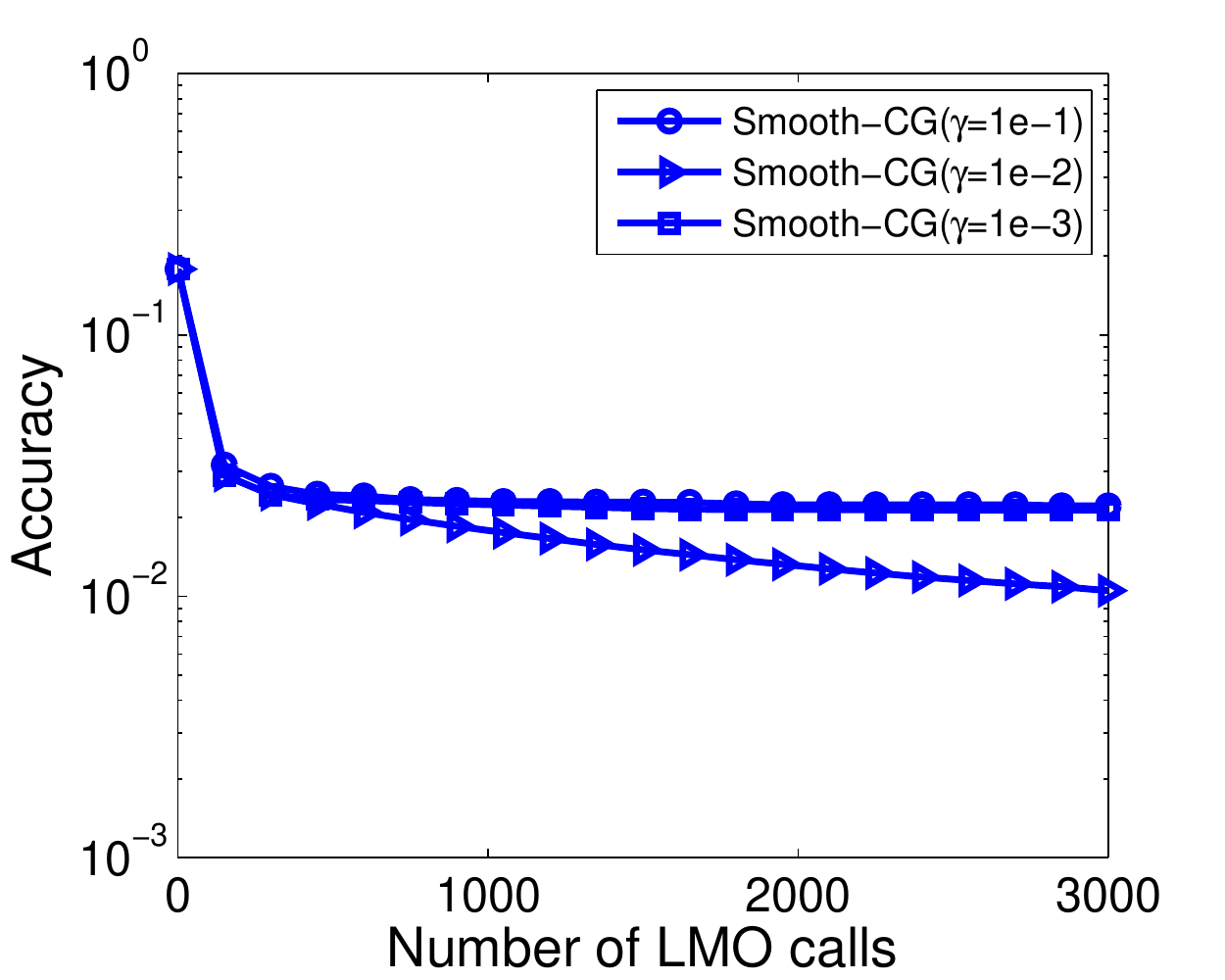}
  \end{minipage}%
   \begin{minipage}[t]{.25\textwidth}
    \includegraphics[scale=0.3]{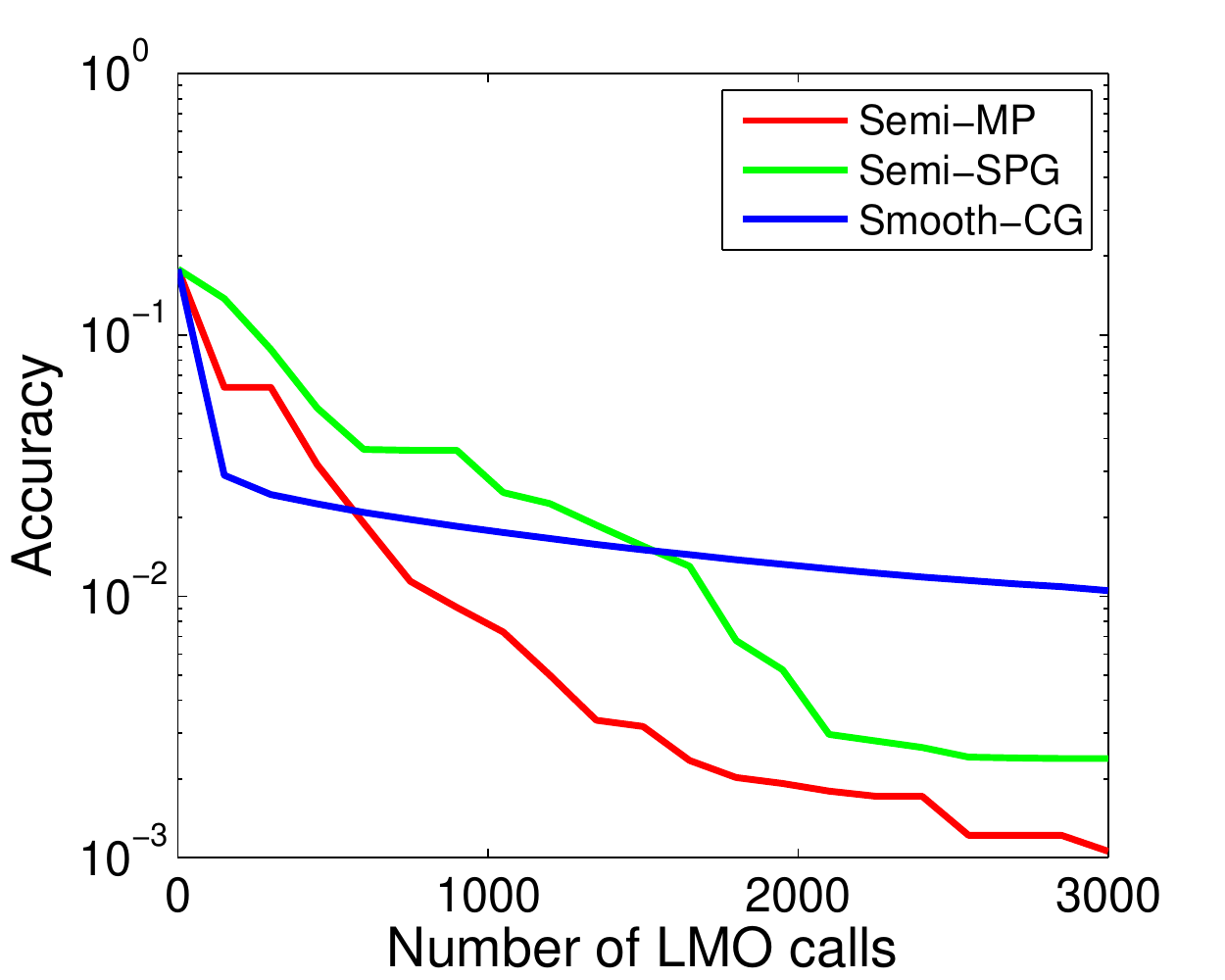}
  \end{minipage}
\caption{Matrix completion on synthetic data($1024\times1024$): optimality gap vs the LMO calls.\\ 
   From left to right:  (a) Semi-MP; (b) Semi-SPG ; (c) Smooth-CG; (d) best of three algorithms.}
   \label{fig:synthetic}
\end{figure}

\subsection{Robust collaborative fitering: $\ell_1$-empirical risk +nuclear norm} 
We consider the collaborative filtering problem, with a nuclear-norm regularisation penalty and an $\ell_1$-empirical risk function:
\begin{equation}\label{matrixcompletionL1}
\min_{x} \frac{1}{|E|}\sum_{(i,j)\in E}|x_{ij}-b_{ij}|+\lambda\|x\|_\nuc.
\end{equation}

\paragraph{Competing algorithms.} We compare the above three candidate algorithm. The smoothed problem for Semi-SPG and Smooth-CG in this case becomes
\begin{equation}\label{smoothL1}
\min _{x,v: \|x\|_\nuc\leq v}f^\gamma(x)+\lambda v, \text{ where } f^{\gamma}(x)=\max_{\|y\|_\infty\leq 1}\left\{ \frac{1}{|E|}\sum_{(i,j)\in E}(x_{ij}-b_{ij})y_{ij}-\frac{\gamma}{2}\|y\|_2^2\right\} .
\end{equation}

Note that in this case, for Smooth-CG, solving the re-optimization problem in (\ref{reoptimize}) at each iteration 
requires computing the full matrix representation for the gradient. For large $t$ and large-scale problems, the computation cost for re-optimization is no longer negligible. However, the Semi-MP and Semi-SPG do not suffer from this limitation since the conditional gradient routines are called for simple quadratic subproblems. For this particular example, we implement the Semi-MP slightly different from the above scheme. We solve the following saddle point reformulation with properly selected $\rho$,
\begin{equation}\label{saddlepointL1}
\min _{\stackrel{x,y, v_1,v_2:}{v_1\geq \|x\|_\nuc,v_2\geq\|y\|_1}}\max_{\|w\|_2\leq 1}  v_2+\lambda v_1+\rho\langle \mathcal{A}x-b-y, w\rangle
\end{equation}
where we use $\mathcal{A}$ to denote the operator $\frac{1}{|E|}P_E$. 
 The semi-structured variational inequality \semivi$(X,F)$ associated with the above saddle point problem is given by 
$X=\{[u=(x,y,w);v=(v_1.v_2)]: \|x\|_\nuc\leq v_1,\|y\|_1\leq v_2,\|w\|_2\leq 1\}$ and $F=[F_u(u);F_v]=[\rho\cA w; -\rho w;\rho(y-\cA x+b);\lambda;1]$. The subdomain $X_1=\{(y,w,v_2):\|y\|_1\leq v_2,\|w\|_2\leq 1\}$ is given by full-prox setup and the subdomain $X_2=\{(x;v_1):\|x\|_\nuc\leq v_1\}$ is given by LMO. By setting both the distance generating functions  as the Euclidean distance, the update of $w$ reduces to the gradient step,  the update of $y$ reduces to the soft-thresholding operator, and the update of $x$ is given by the composite conditonal gradient routine. In our experiment, the factor $\rho$ is updated adaptively in such a way that the back-projection step does not increase the objective function value. We set the stepsizes $\gamma_t$ along the iterations using line-search. All in all, the \spmp algorithm (Semi-MP) is fully automatic, and does not require tuning of any parameter.

We run the above three algorithms on the the small and medium MovieLens datasets. The small-size dataset consists of 943 users and 1682 movies with about 100K ratings,while the medium-size dataset consists of 3952 users and 6040 movies with about 1M ratings. We follow~\cite{pierucci2014smoothing} to set the regularisation parameters. We randomly pick 80\% of the entries to build the training dataset, and compute the normalized mean absolute error (NMAE) on the remaining test dataset.  For Smooth-CG, we carry out the algorithm with different smoothing parameters, ranging from $\{1e-3, 1e-2,1e-1, 1e0\}$ and select the one with the best performance. For the Semi-SPG algorithm, we adopt the best smoothing parameter found in Smooth-CG.  We use two different strategies to control the number of LMO calls at each iteration, i.e. the accuracy of the proximal mapping for both Semi-SPG and Semi-MP, which are a) fixed inner CG steps and b) decaying $\epsilon_t=c/t$ as the theory suggested.  We report in Fig.~\ref{fig:small} and Fig.~\ref{fig:medium} the performance of each algorithm under different choice of parameters and the overall comparison of objective value and NMAE on test data in Fig.~\ref{fig:robustCF}. 

\begin{figure} [!ht]
  \begin{minipage}[t]{.25\textwidth}
    \includegraphics[scale=0.3]{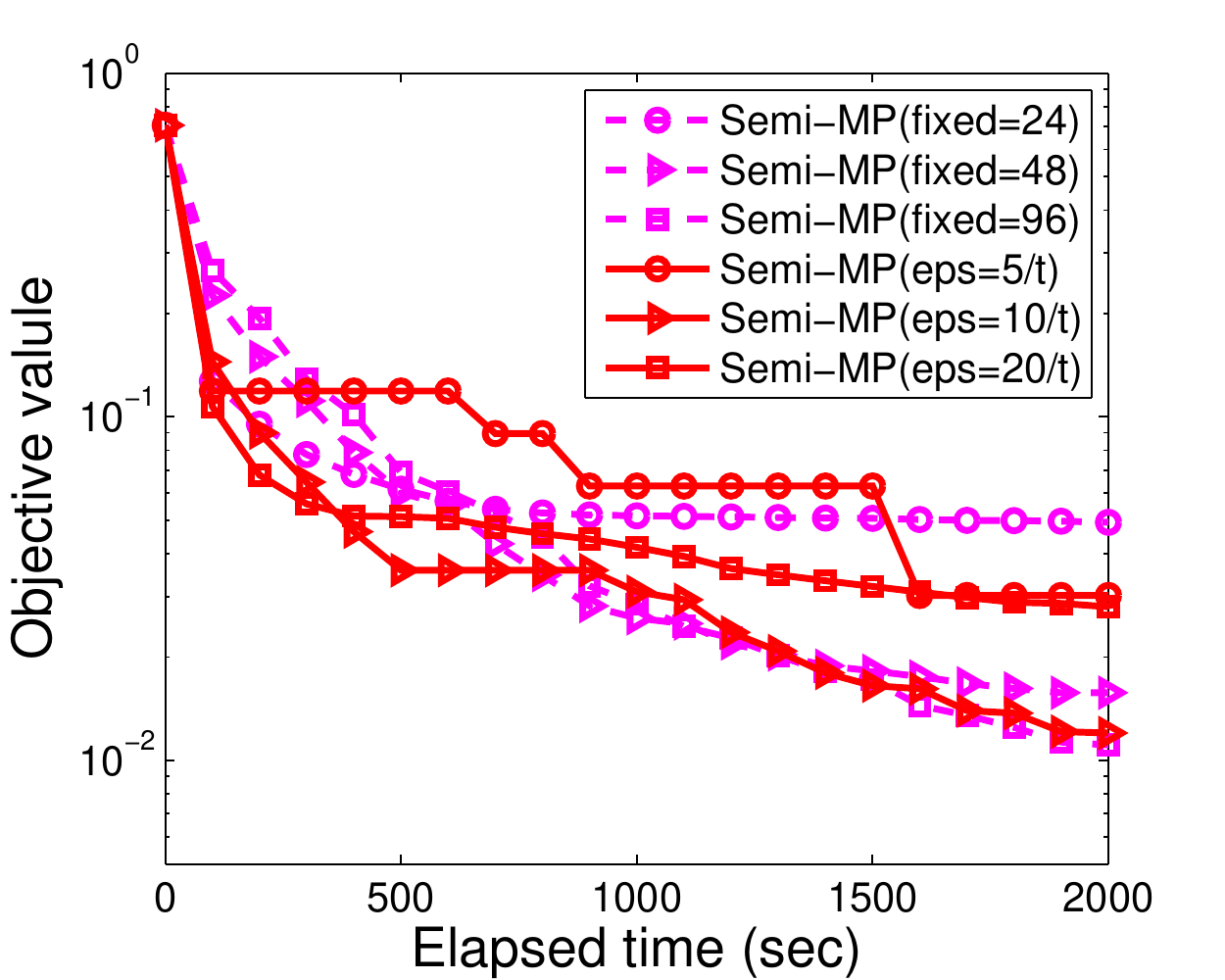}
  \end{minipage}%
  \begin{minipage}[t]{.25\textwidth}
    \includegraphics[scale=0.3]{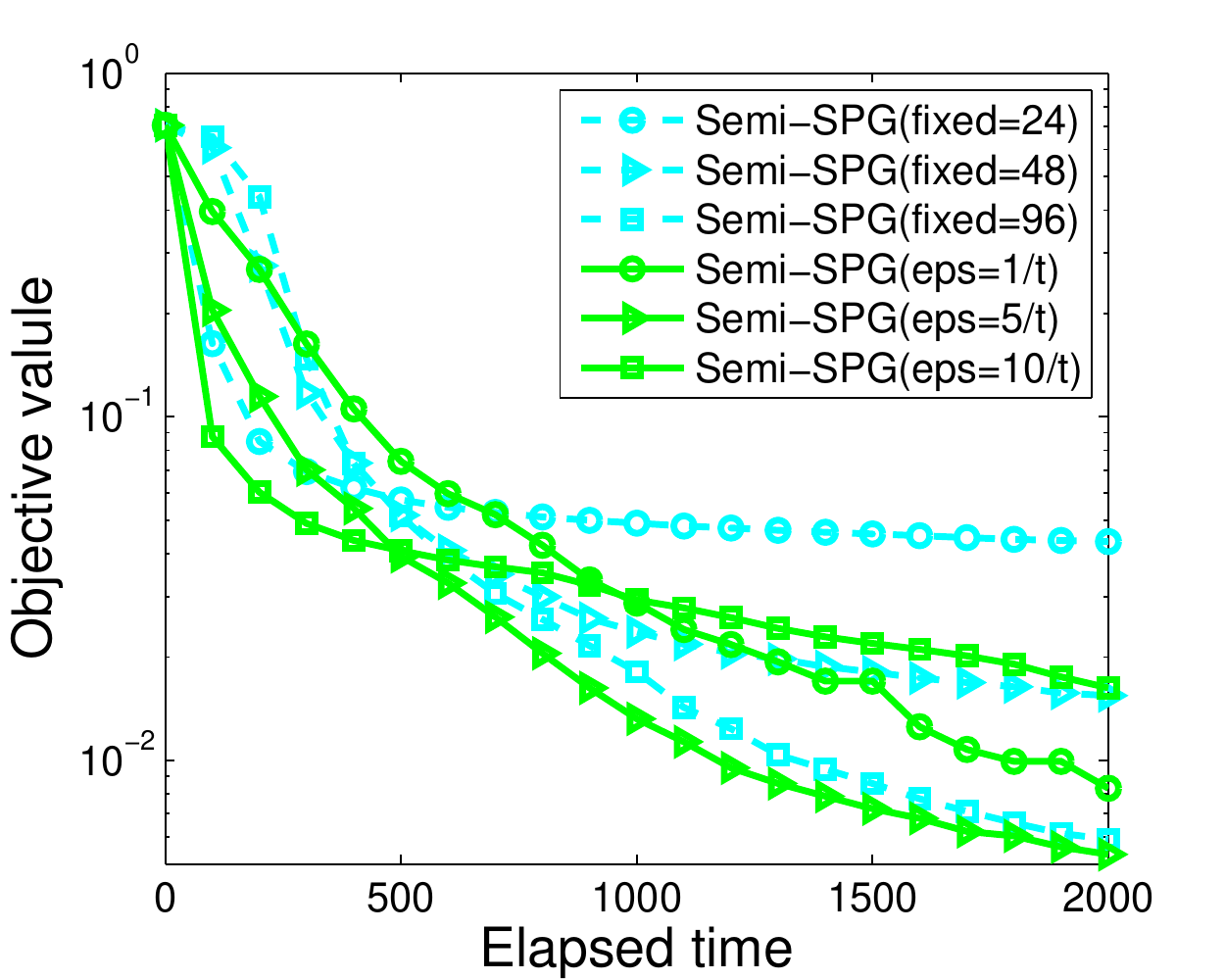}
  \end{minipage}%
  \begin{minipage}[t]{.25\textwidth}
    \includegraphics[scale=0.3]{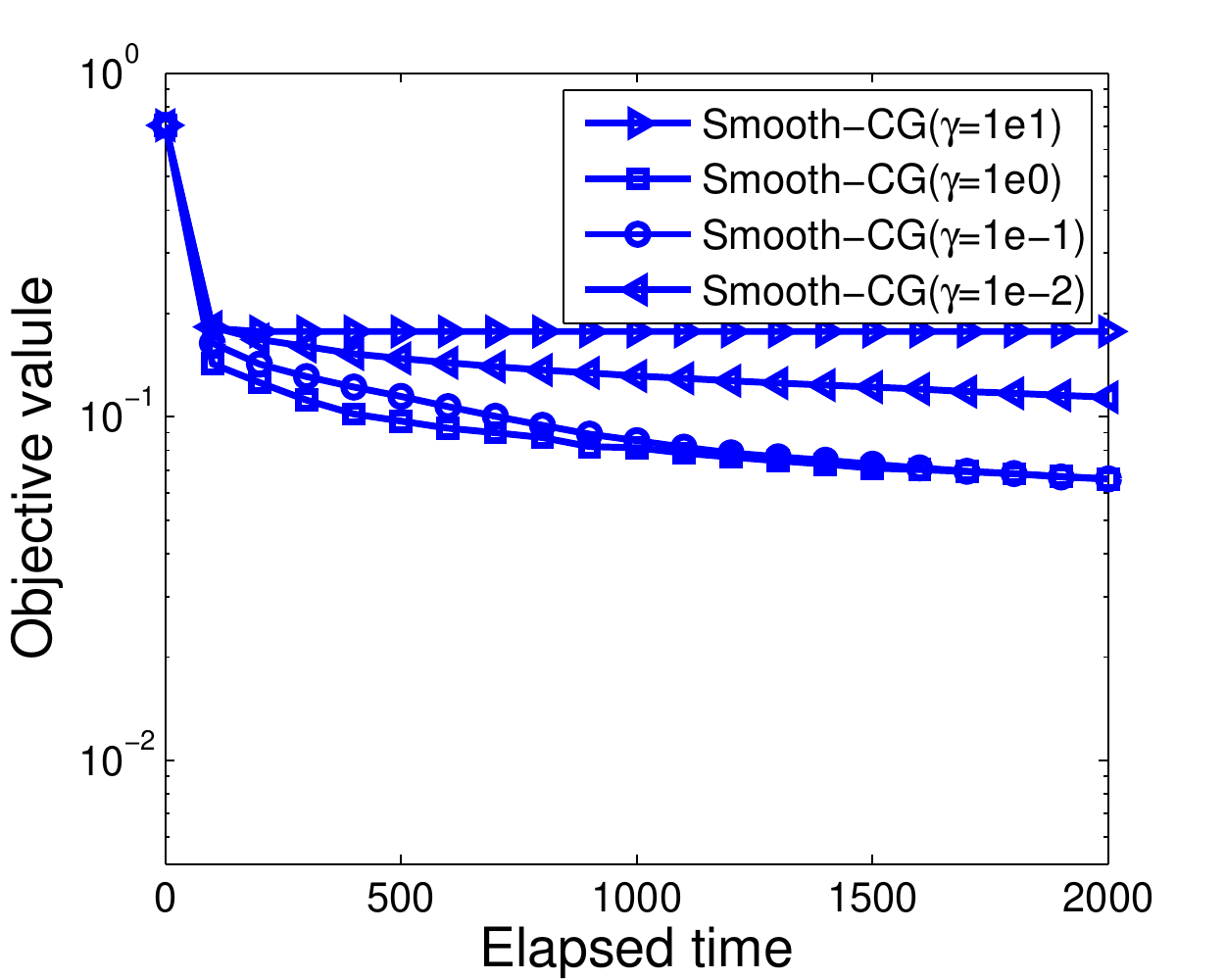}
  \end{minipage}%
  \begin{minipage}[t]{.25\textwidth}
    \includegraphics[scale=0.3]{Small_100K_L1_obj.pdf}
  \end{minipage}%
\caption{Robust collaborative filtering on MovieLens 100K: objective function vs elapsed time. \\
   From left to right:  (a) Semi-MP; (b) Semi-SPG ; (c) Smooth-CG; (d) best of three algorithms.} 
\label{fig:small}
\end{figure}

\begin{figure} [!ht]
  \begin{minipage}[t]{.25\textwidth}
    \includegraphics[scale=0.3]{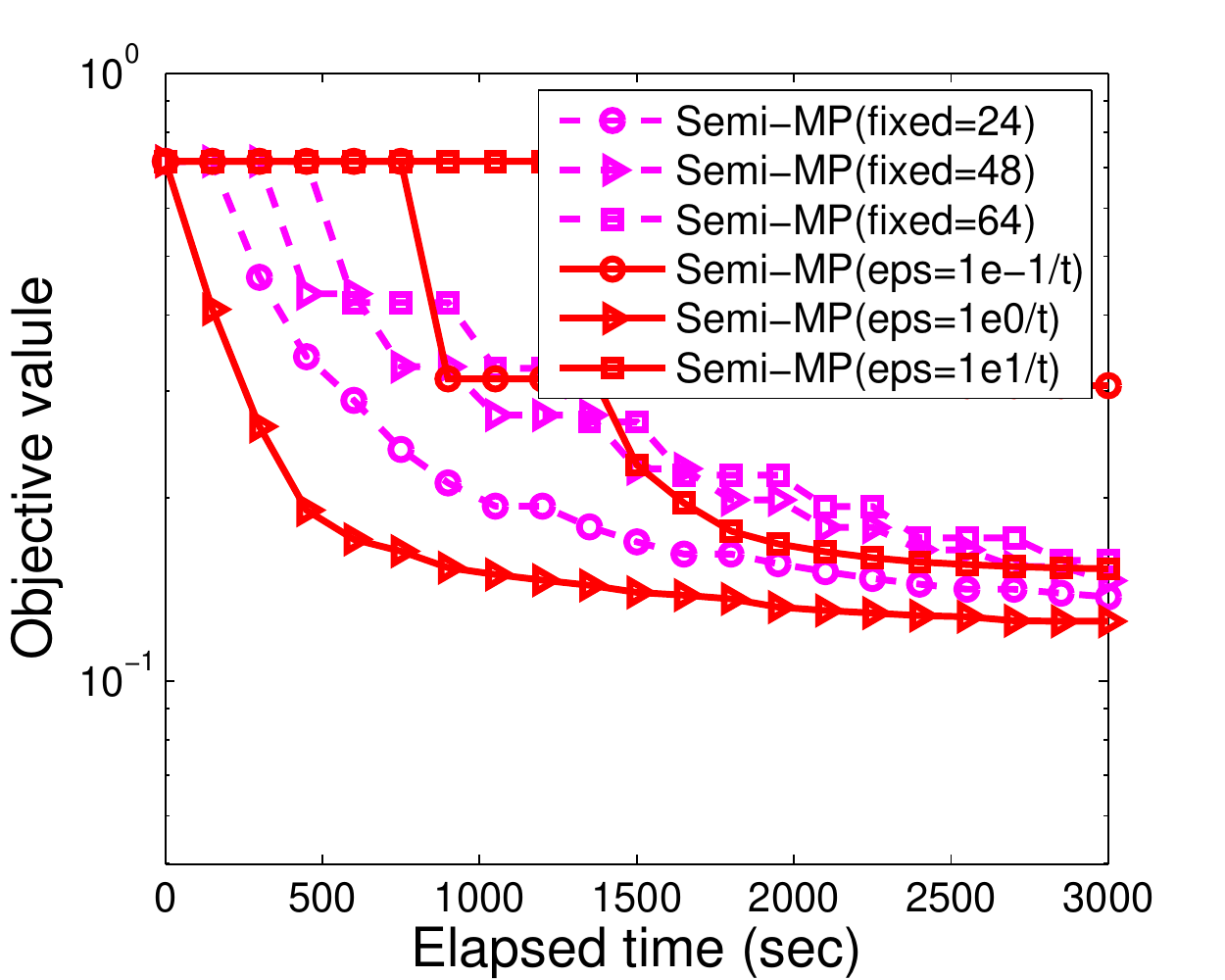}
  \end{minipage}%
  \begin{minipage}[t]{.25\textwidth}
    \includegraphics[scale=0.3]{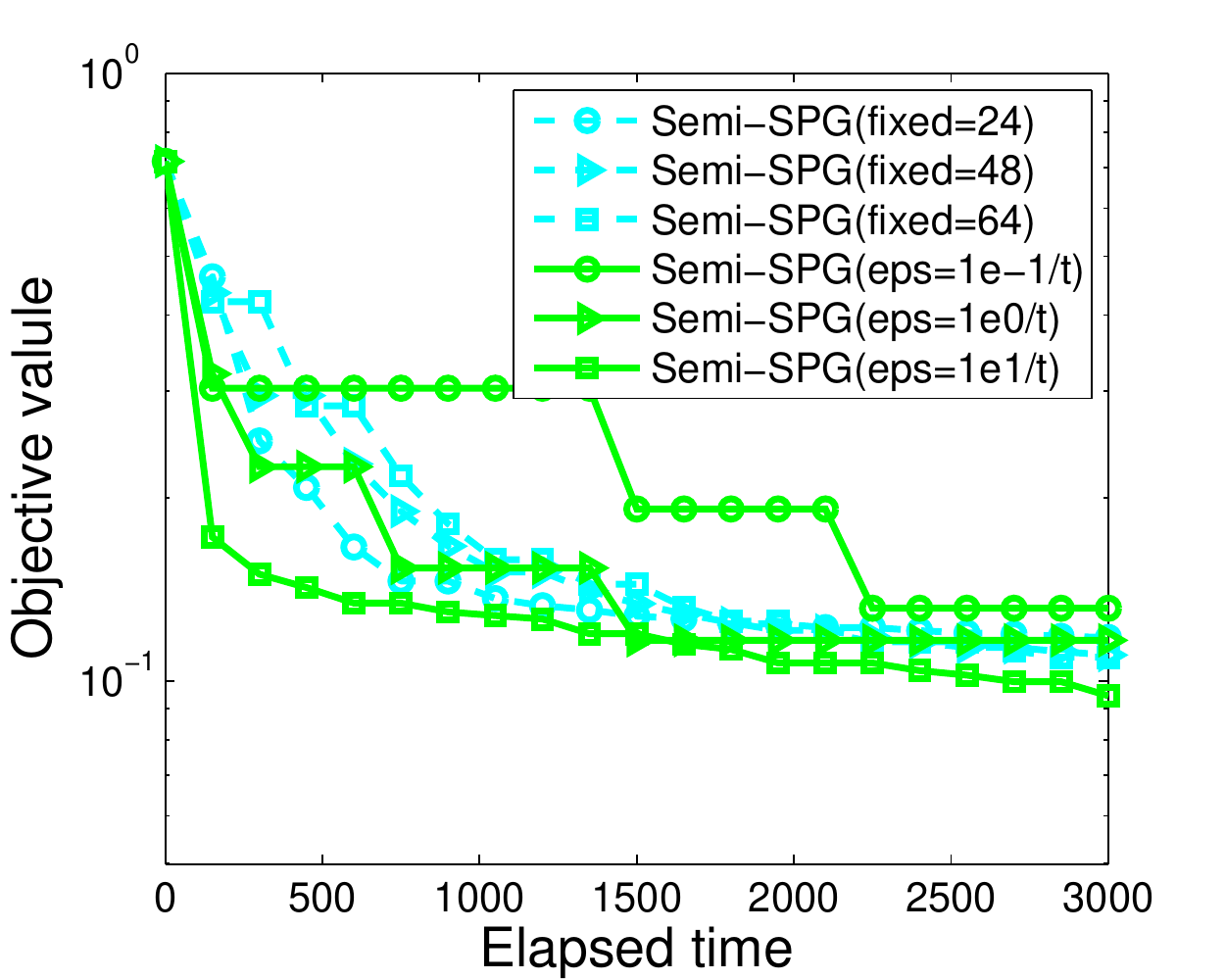}
  \end{minipage}%
  \begin{minipage}[t]{.25\textwidth}
    \includegraphics[scale=0.3]{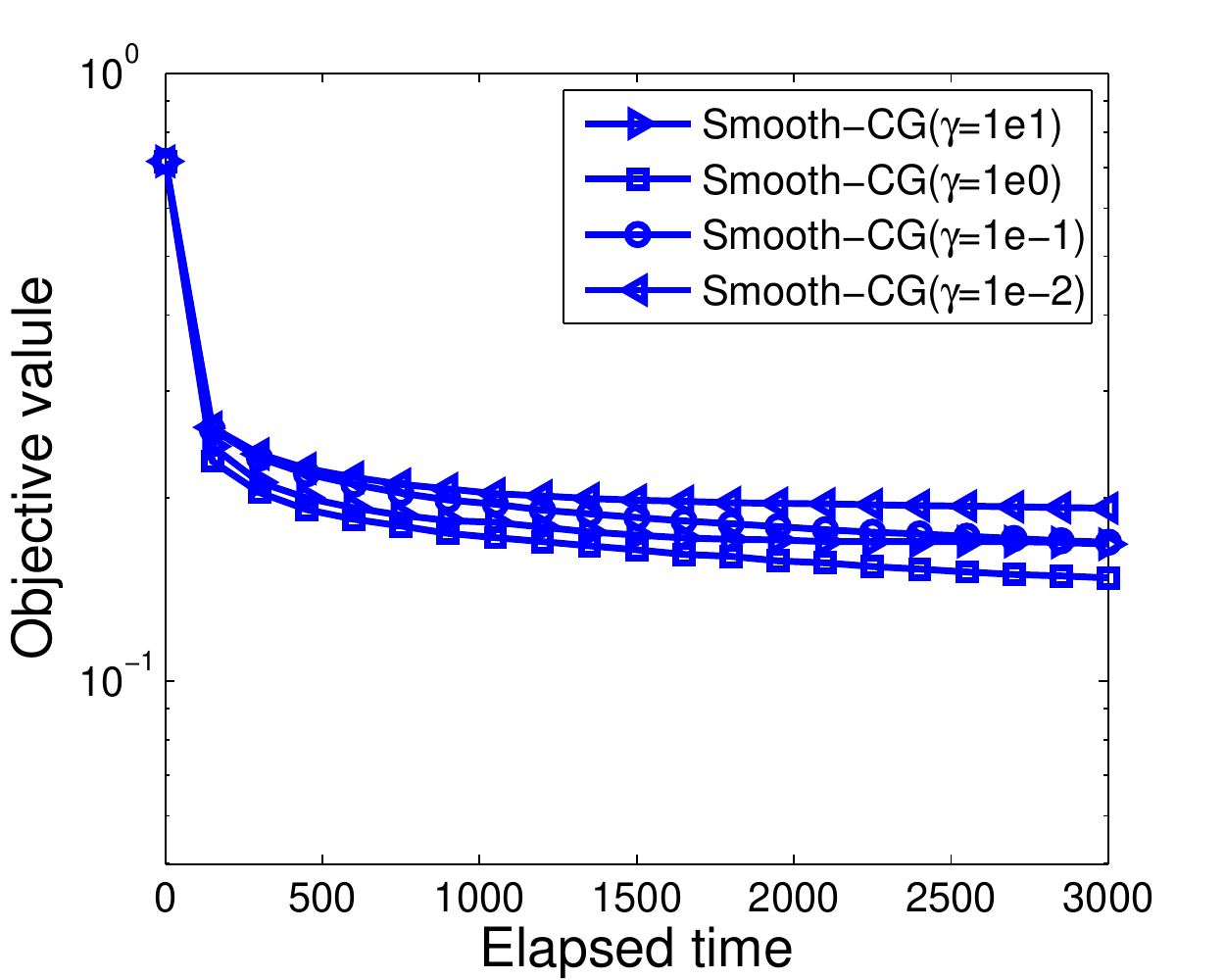}
  \end{minipage}%
  \begin{minipage}[t]{.25\textwidth}
    \includegraphics[scale=0.3]{Medium_1M_L1_obj.pdf}
  \end{minipage}%
\caption{Robust collaborative filtering on MovieLens 1M: objective function vs elasped time. \\
   From left to right:  (a) Semi-MP; (b) Semi-SPG ; (c) Smooth-CG; (d) best of three algorithms.} 
\label{fig:medium}
\end{figure}

\begin{figure}[!ht]
  \begin{minipage}[t]{.25\textwidth}
    \includegraphics[scale=0.28]{Small_100K_L1_obj.pdf}
  \end{minipage}%
  \begin{minipage}[t]{.25\textwidth}
    \includegraphics[scale=0.28]{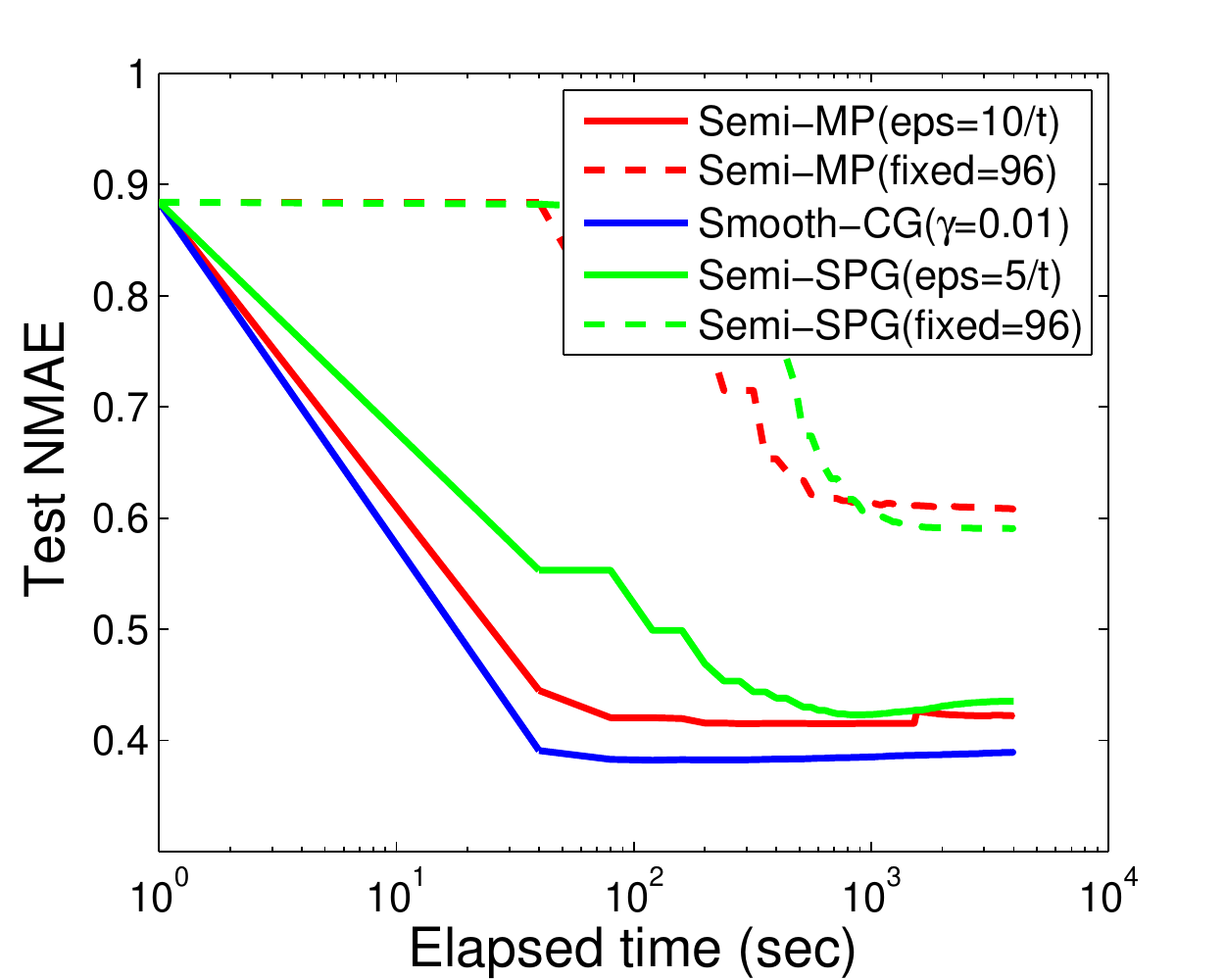}
  \end{minipage}%
  \begin{minipage}[t]{.25\textwidth}
    \includegraphics[scale=0.28]{Medium_1M_L1_obj.pdf}
  \end{minipage}%
 \begin{minipage}[t]{.25\textwidth}
    \includegraphics[scale=0.28]{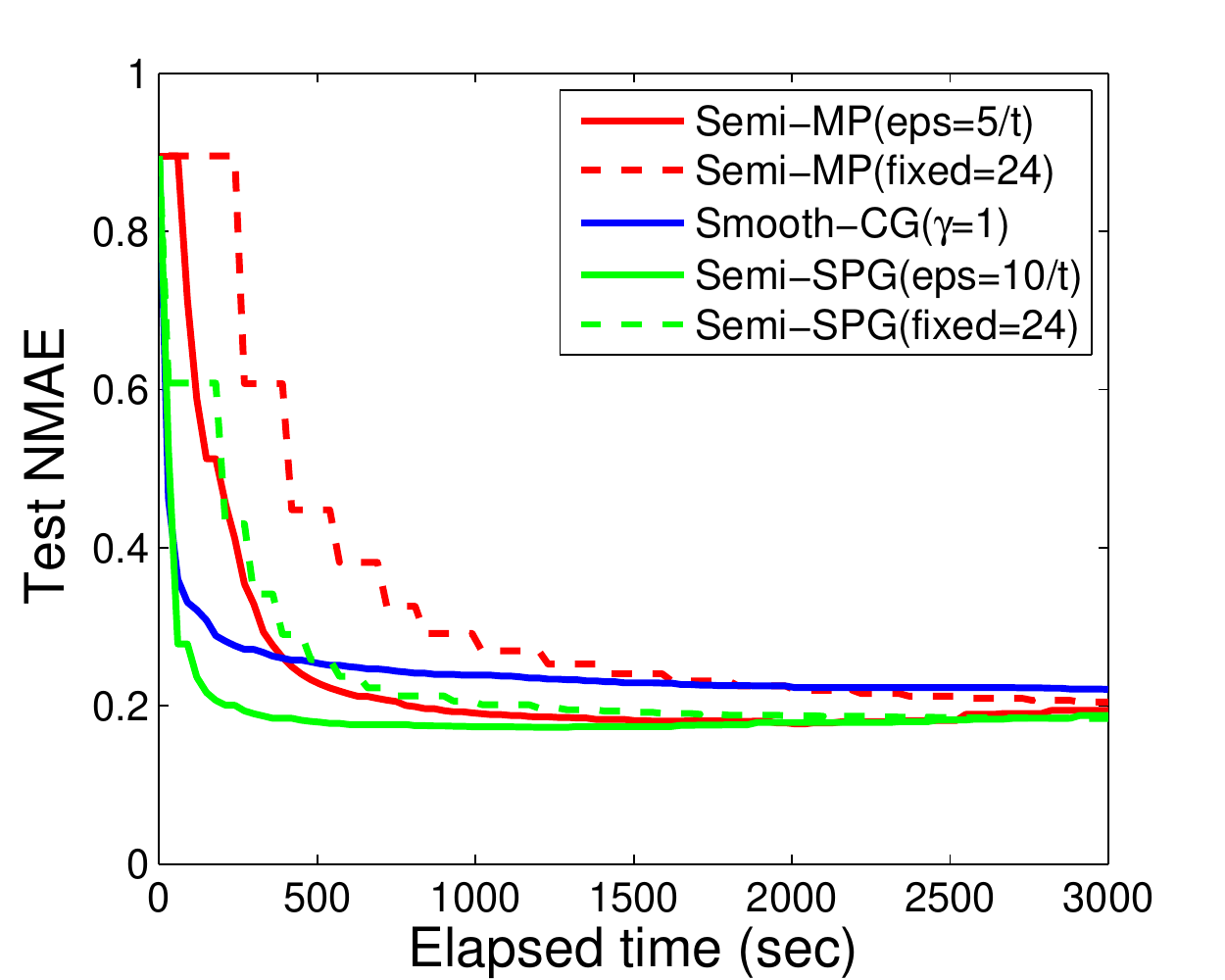}
  \end{minipage}
\caption{Robust collaborative filtering on Movie Lens: objective function and test NMAE against elapsed time.   From left to right: (a) MovieLens 100K objective; (b) MovieLens 100K test NMAE; (c) MovieLens 1M objective; (d) MovieLens 1M test NMAE.}
\label{fig:robustCF}
\end{figure}

 In Fig.~\ref{fig:small} and Fig.~\ref{fig:medium}, we can see that using fixed inner CG steps sometimes achieve comparable performance as using the decaying epsilon $\epsilon_t$. In Fig.~\ref{fig:robustCF}, we can see that Semi-MP clearly outperforms Smooth-CG, while it is competitive with Semi-SPG. In the large-scale setting, Semi-MP achieves better objective as well as test NMAE compared to Smooth-CG.

\subsection{Link prediction: hinge loss + $\ell_1$-norm + nuclear norm}
We consider the following model for the link prediction problem, 
\begin{equation}\label{hingeloss}
\min _{x\in\bR^{m\times n}} \frac{1}{|E|}\sum_{(i,j)\in E}\max\left(1-(b_{ij}-0.5)x_{ij},0\right)+\lambda_1\|x\|_1+\lambda_2\|x\|_\nuc
\end{equation}

This example is more complicated than the previous two examples since it has not only one nonsmooth loss function but also two regularization terms. Applying the smoothing-CG or Semi-SPG would require to build two smooth approximations, one for hinge loss term and one for $\ell_1$ norm term. Therefore, we consider another alternative approach, Semi-LPADMM, where we apply the linearized preconditioned ADMM algorithm by solving proximal mapping through conditional gradient routines. Up to our knowledge, ADMM with early stopping is not well-analyzed in literature, but intuitively as long as the accumulated error is controlled sufficiently, the variant will converge. 

We conduct experiments on a binary social graph data set called Wikivote, which consists of 7118 nodes and 103,747 edges. Since the computation cost of these two algorithms mainly come from the LMO calls, we present in below the performance in terms of number of LMO calls. For the first set of experiments, we select top 1024 highest degree users from Wikivote and run the two algorithms on this small dataset with different strategies for the inner LMO calls. 

In Fig.~\ref{fig:linkpredict}, we observe that the Semi-MP is less sensitive to the inner accuracies of prox-mappings compared to the ADMM variant, which sometimes stop progressing if the prox mapping of early iterations are not solved with sufficient accuracy.  Another observation is that in this example, the second strategy, which essentially saves the use of LMOs, works better in the long run than using fixed number of LMOs. The results indicate again on the full dataset  again indicates that our algorithm performs better than the semi-proximal variant of ADMM algorithm.

\begin{figure}[!ht]
  \begin{minipage}[t]{.33\textwidth}
    \includegraphics[scale=0.38]{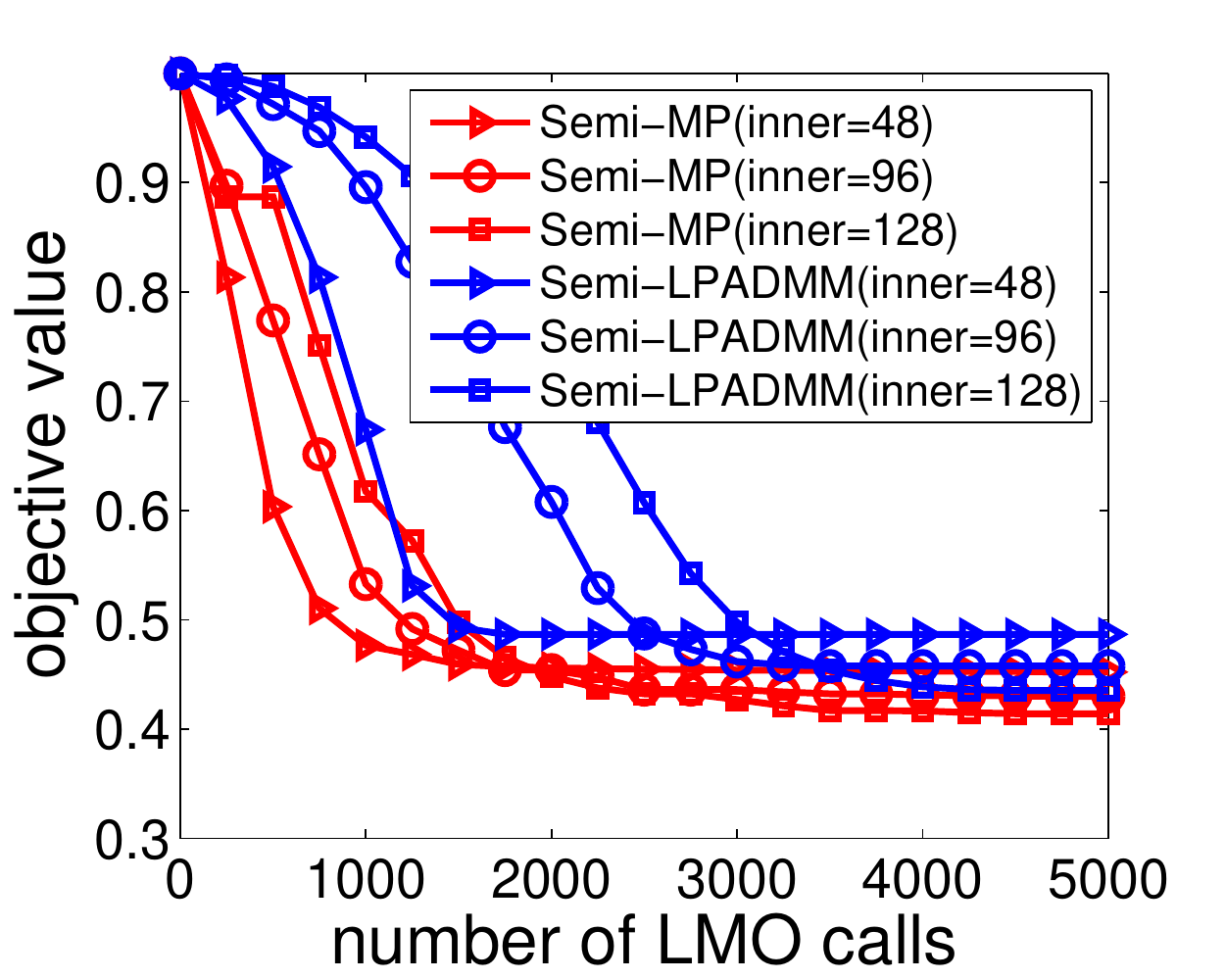}
  \end{minipage}%
  \begin{minipage}[t]{.33\textwidth}
    \includegraphics[scale=0.38]{decay.pdf}
  \end{minipage}%
 \begin{minipage}[t]{.33\textwidth}
    \includegraphics[scale=0.38]{wikivote_two.pdf}
  \end{minipage}
\caption{Link prediction on Wikivote: objective function value against the LMO calls. 
From left to right: (a)Wikivote(1024) with fixed inner steps; (b) Wikivote(1024) with $\epsilon_t=c/t$;  (c) Wikivote(full)}
\label{fig:linkpredict}
\end{figure}
%

\end{document}